\documentclass[jmp,preprint]{revtex4-1}

\draft 

\usepackage{tikz}
\usetikzlibrary{decorations.pathreplacing,angles,quotes}

\usetikzlibrary{matrix,arrows,decorations.pathmorphing}
\usepackage{tikz-cd}
\usetikzlibrary{arrows}
\usepackage{verbatim}

\usepackage{mathrsfs}

\usepackage{epsfig}
\usepackage{amsmath,amsfonts,amsthm,amscd, amssymb}
\usepackage{latexsym}
\usepackage{enumerate}

\newcommand*{\rom}[1]{\expandafter\@slowromancap\romannumeral #1@}
\newcommand{\Proof}[1]{\begin{proof} #1 \end{proof}}
\newtheorem{thm}{Theorem}[section]
\newtheorem*{theorem*}{Theorem}
\newtheorem*{question*}{Question}

\newtheorem{pro}[thm]{Proposition}
\newcommand{\Pro}[1]{\begin{pro} #1\end{pro}}
\newtheorem{cor}[thm]{Corollary}
\newcommand{\Cor}[1]{\begin{cor} #1 \end{cor}}
\newtheorem{lem}[thm]{Lemma}
\newcommand{\Lem}[1]{\begin{lem} #1 \end{lem}}

\newtheorem{construction}[thm]{Construction}
\newtheorem{defn}[thm]{Definition}
\newcommand{\Def}[1]{\begin{defn}{\rm #1}\end{defn}}
\newtheorem{rem}[thm]{Remark}
\newcommand{\Rem}[1]{\begin{rem} #1 \end{rem}}
\newtheorem{ex}[thm]{Example}
\newcommand{\Ex}[1]{\begin{ex} #1 \end{ex}}

\newtheorem{obs}[thm]{Observation}

\newtheorem{note}[thm]{Note}

\newtheorem{nota}[thm]{Notation}

\newcommand{\catTop}{\mathbf{Top}}
\newcommand{\catFin}{\mathbf{Fin}}

 \newcommand{\Spin}{\text{Spin}}

\newcommand{\tightoverset}[2]{%
\mathop{#2}\limits^{\vbox to -.5ex{\kern-0.25ex\hbox{$#1$}\vss}}}

\newcommand{\Hom}{\text{Hom}}

\newcommand{\Sol}{\text{Sol}}

\newcommand{\sym}{\text{sym}}

\newcommand{\idy}{\text{id}}

\newcommand{\Det}{\text{det}}

\newcommand{\Rep}{\text{Rep}}

\newcommand{\sq}{\text{sq}} 
\newcommand{\st}{\text{st}}

\newcommand{\cl}{\text{cl}}

\newcommand{\hquo}{\ensuremath{/ \! /}} 

\newcommand{\twobytwo}[4]{\ensuremath{\left( \begin{matrix}
#1 & #2 \\
#3 & #4
\end{matrix}  \right) }}

\newcommand{\set}[1]{\ensuremath{ \lbrace #1 \rbrace }}

\newcommand{\Span}[1]{\ensuremath{ \langle #1 \rangle }}

\newcommand{\com}{\text{com}}  

\newcommand{\Gr}{\text{Gr}}

\newcommand{\dt}{\bullet}

\newcommand{\SO}{\text{SO}}

\newcommand{\bi}{\mathfrak{b}}

\newcommand{\ZZ}{\mathbb{Z}}

\newcommand{\RR}{\mathbb{R}}
\newcommand{\CC}{\mathbb{C}}
\newcommand{\NN}{\mathbb{N}}

\newcommand{\sS}{\mathbb{S}}

\newcommand{\hH}{\mathfrak{H}}

\usepackage[parfill]{parskip}
\numberwithin{equation}{section}
\renewcommand{\theequation}{\arabic{section}.\arabic{equation}}
\usepackage{multirow}
\usepackage{tabularx}

\newcommand\rev[1]{{\color{black}#1}}
\newcommand\revtwo[1]{{\color{black}#1}}
\newcommand\revthr[1]{{\color{black}#1}}

\newtheorem{innercustomthm}{Theorem}
\newenvironment{customthm}[1]
  {\renewcommand\theinnercustomthm{#1}\innercustomthm}
  {\endinnercustomthm}

\newtheorem{innercustomcor}{Corollary}
\newenvironment{customcor}[1]
  {\renewcommand\theinnercustomcor{#1}\innercustomcor}
  {\endinnercustomcor}

\begin{document}


\title{Commutative $d$-torsion $K$-theory and its applications} 



\author{Cihan Okay}
\email[]{cihan.okay@bilkent.edu.tr}
\affiliation{Bilkent University, Ankara, Turkey}


\date{\today}

\begin{abstract}
Commutative $d$-torsion $K$-theory is a variant of topological $K$-theory constructed from commuting unitary matrices of order dividing $d$.
Such matrices appear as solutions of linear constraint systems that play a role in the study of 
quantum contextuality and in applications to operator-theoretic problems motivated by quantum information theory.
Using methods from stable homotopy theory we modify commutative $d$-torsion $K$-theory
into a cohomology theory which can be used for studying operator solutions of linear constraint systems.  This provides an interesting connection between stable homotopy theory and  quantum  information theory.
\end{abstract}


\maketitle 


\section{Introduction}

Commuting unitary matrices can be assembled into a generalized cohomology theory called commutative $K$-theory, a variant of topological $K$-theory first introduced in \cite{AGLT}. 
This theory can be further modified by restricting to matrices whose order divides $d$. The  resulting  cohomology theory  will be referred to as {\it commutative $d$-torsion $K$-theory}.  Such matrices also play a significant role in quantum theory, especially in foundational areas concerning quantum contextuality \cite{kochen,bell1966problem} and linear constraint systems in the study of non-local games \cite{cleve2014characterization}. The goal of this paper is to make this connection precise.
We 
\revthr{construct} another generalized cohomology theory\revthr{, called the $C(d,m)$-cohomology,}   obtained from commutative $d$-torsion $K$-theory, which is tailored for studying operator solutions of linear constraint systems.   We expect that stable homotopical methods introduced in this paper will provide further insight into operator-theoretic problems motivated by quantum information theory.

The cohomology theories studied in this paper are based on a classifying space construction   introduced in \cite{ACT12}. 
We write $B(\ZZ/d,G)$ to denote the classifying space of a topological group $G$ constructed from tuples of pairwise commuting group elements where each group element has order dividing $d$, i.e., pairwise commuting $d$-torsion group elements. When $G$ is the unitary group $U(m)$   this classifying space is constructed from tuples $(A_1,A_2,\cdots,A_n)$ of matrices satisfying 
$$
  A_iA_j=A_jA_i \;\; \text{ and }\;\;    (A_i)^d=I_m
$$
where $I_m$ is the $m\times m$ identity matrix.
Such matrices also appear as solutions to a linear constraint system    specified by an equation $Mx=b$ where $M$ is an $r\times c$ matrix over the additive group $\ZZ/d$ of integers modulo $d$. 
An operator solution consists of $d$-torsion $m\times m$ unitary matrices $A_1,A_2,\cdots ,A_c$  that satisfy  
$$
A_1^{M_{k1}} A_2^{M_{k2}} \cdots A_c^{M_{kc}} = e^{2\pi ib_k /d} I_m\;\;\;\text{ for all }\;\; 1\leq k\leq r,
$$
and $A_iA_j=A_jA_i$ whenever $M_{ki}$ and $M_{kj}$ are both non-zero. The data of a linear constraint system can be packaged as a pair $(\hH,\tau)$, where $\hH$ is a hypergraph with a vertex set $V=\set{v_1,v_2,\cdots,v_c}$ and an edge set $E=\set{e_1,e_2,\cdots,e_r}$. Here $\tau$ is the function $E\to \ZZ/d$ defined by $\tau(e_k)=b_k$.  An operator solution $\set{A_i}$ can be regarded as a function $T:V\to U(m)$ where $T(v_i)=A_i$.
The homotopical approach initiated in \cite{Coho,okay2020homotopical} associates a $2$-dimensional CW complex $X$, called a topological realization, to the hypergraph $\hH$ and the function $\tau$ represents a $2$-dimensional cohomology class $[\tau]\in H^2(X,\ZZ/d)$. 

In this paper, we refine this approach by interpreting an operator solution as a map of topological spaces. For this a quotient space $\bar B(\ZZ/d,G)$ of the classifying space $B(\ZZ/d,G)$ is introduced. An operator solution over $G=U(m)$ can be turned into a map, defined up to homotopy,
$$
f_T:X \to \bar B(\ZZ/d,G).
$$
Although,   our motivation comes from an urge to understand operator solutions of linear constraint systems, the classifying space $B(\ZZ/d,G)$ and its variants are of independent interest to algebraic topologists; see for instance  \cite{O16,AG15,cohen2016survey,AGV17,OW18,DV18,ramras2018homological,okay2020commutative}.

A generalized cohomology theory is represented by a spectrum. Following \cite{gritschacher2019commuting} we show that $B(\ZZ/d,U)$, where $U$ is the stable unitary group, is an infinite loop space and thus specifies a spectrum.  
This spectrum  turns out to be stably equivalent to $ku\wedge B\mu_d$ (Proposition \ref{pro:InfLoop}), where $\mu_d=\set{e^{2\pi ik /d}|\;1\leq k\leq d}$ and $ku$ is the connective complex $K$-theory spectrum.  Commutative $d$-torsion $K$-theory is the generalized cohomology theory associated to this spectrum. Both the spectrum and the associated cohomology theory will be denoted by $k\mu_d$. 

For applications to linear constraint systems we introduce a stabilized version of the quotient space $\bar B(\ZZ/d,U(m))$. The usual stabilization process cannot be carried out in a straightforward manner. However, by working in the homotopy category of spectra we 
\revthr{construct} a spectrum\revthr{, denoted by $C(d,m)$,  obtained from $k\mu_d$. It turns out that the resulting cohomology theory is represented by a product of two Eilenberg--MacLane spaces, thus  can be constructed more directly.  
} 
\revthr{The} infinite loop space $\bar B(d,m)$ \revthr{associated to the spectrum $C(d,m)$} admits a map
$$
\bar\iota_m:\bar B(\ZZ/d,U(m))\to \bar B(d,m).
$$
This space   comes with a canonical cohomology class $\gamma_m^\sS$ in $H^2(\bar B(d,m),\ZZ/d)$.  
By construction homotopy groups of $C(d,m)$ are concentrated in dimensions $i=1,2$ and we show that there is an exact sequence
$$
0\to \pi_2 C(d,m) \to \ZZ/d \xrightarrow{\times m} \ZZ/d \to \pi_1 C(d,m) \to 0.
$$ 
The kernel consists of the subgroup  $(\ZZ/d)_m$ of $m$-torsion elements.
Using  the Atiyah--Hirzebruch spectral sequence we describe $C(d,m)$-cohomology of a space:

\begin{customthm}{\ref{thm:Thm-stable}}
Let $X$ be a connected 
CW complex. There is a commutative diagram
\begin{equation}\label{eq:diag-thm}
\begin{tikzcd}
& H^2(X,(\ZZ/d)_m) \arrow[d,hook]\arrow[r,"(i_m)_*"] &  H^2(X,\ZZ/d) \arrow[d,equal] \\
k\mu_d(X) \arrow[d] \arrow{r}{\zeta} & C(d,m)(X) \arrow[d,two heads] \arrow{r}{\cl} & H^2(X,\ZZ/d) \\ 
H^1(X,\ZZ/d) \arrow[r,"(\pi_m)_*"] & H^1(X,\frac{\ZZ/d}{m\ZZ/d}) &
\end{tikzcd}
\end{equation}
where  $\cl(f)=f^*(\gamma_m^\sS)$, the image of $\zeta$ is contained in the kernel of $\cl$, and the middle column is an exact sequence 
\revtwo{which (non-canonically) splits as follows
$$
C(d,m)(X) \cong H^1(X,\frac{\ZZ/d}{m\ZZ/d}) \oplus H^2(X,(\ZZ/d)_m) .
$$}
\end{customthm}

Going back to linear constraint systems we show that the  $C(d,m)$-cohomology informs us about the properties of operator solutions over $U(m)$. To an operator solution we associate the class $[f]$ of the composite map 
$$
f:X \xrightarrow{f_T} \bar B(\ZZ/d,U(m)) \xrightarrow{\bar \iota_m} \bar B(d,m)
$$
in the $C(d,m)$-cohomology of $X$. It turns out that $\cl(f)=0$ if and only if the linear constraint system has a solution over $U(1)$, also known as a scalar solution.

\begin{customcor}{\ref{cor:stable-LCS}}
Let $(\hH,\tau)$ be a linear constraint system  over $\ZZ/d$ and $X$ be a topological realization for $\hH$.
\begin{enumerate}
\item If \rev{$(\hH,\tau)$ has an operator solution and} $H^2(X,(\ZZ/d)_m)=0$ then  $(\hH,\tau)$   has a scalar solution.
\item If $d$ and $m$ are coprime then $C(d,m)(X)=0$. In particular,  $(\hH,\tau)$    has a scalar solution if it has an operator solution over $U(m)$.
\item If $\pi_1(X)$ is trivial and $[\tau]\not= 0$ then $(\hH,\tau)$ does not have an operator solution over $U(m)$ for any \rev{$m\geq 2$}. 

\end{enumerate}
\end{customcor}

The most famous example of a linear constraint system,     which does not admit a scalar solution, is the Mermin square construction   \cite{Mermin}. This linear constraint system, defined over $\ZZ/2$, admits an operator solution in $U(2^n)$ for $n\geq 2$. A topological realization for the Mermin square linear constraint system can be chosen to be a torus $X=S^1\times S^1$ with a certain cell structure (Figure \ref{fig:Msq}). 
Then an operator solution   specifies a class in the $C(2,2^n)$-cohomology of the torus
$$
M_n \in C(2,2^n)(S^1\times S^1).
$$
We refer to this class as the Mermin class. In addition, we show that the Mermin class can also be   identified with the generator of  $\pi_2C(2,2^n)=\ZZ/2$. 
There is also a real version of these constructions which works for the orthogonal group $O(m)$. In this case certain generalized cohomology classes can be realized as symmetry-protected topological phases  (\S \ref{sec:ex-real}).

The paper is organized as follows. In \S \ref{sec:ComdTorK} we introduce the classifying space $B(\ZZ/d,G)$ and the type of principal bundles  classified by this space. $\Gamma$-spaces are used to describe the spectrum $k\mu_d$ and Proposition \ref{pro:InfLoop} informs us about its stable homotopy type. Low dimensional homotopy groups are described in \S \ref{sec:lowdim}. 
The quotient space $\bar B(\ZZ/d,G)$ and the spectrum $C(d,m)$ are introduced in \S \ref{sec:Cdm}. In this section we prove Theorem \ref{thm:Thm-stable}, which describes the $C(d,m)$-cohomology of a space. Applications of $C(d,m)$-cohomology are discussed in \S \ref{sec:LCS}, where we introduce linear constraint systems and a topological interpretation of operator solutions. Proposition \ref{pro:class-G} provides a computation of pointed homotopy classes of maps $X\to \bar B(\ZZ/d,G)$ when $X$ is a $2$-dimensional CW complex. Applications to linear constraint systems are given in Corollary \ref{cor:stable-LCS}. The Mermin class is constructed in this section. \rev{In Appendix \ref{app:real-case} we introduce the real versions of these spectra obtained from the orthogonal group.}

\section{Commutative $d$-torsion $K$-theory}\label{sec:ComdTorK}

In this section we introduce a new generalized cohomology theory obtained as a variant of commutative $K$-theory introduced in \rev{\cite{AG15,AGLT}}. Commutative $K$-theory has nice properties such as the spectrum $ku_\com$ representing the theory is stably equivalent to $ku\wedge \CC P^\infty$ as proved in \cite{Gri17}, where $ku$ is the connective complex $K$-theory spectrum. For the $d$-torsion case the spectrum representing the cohomology theory is denoted by $k\mu_d$. It is constructed from commuting unitary matrices whose eigenvalues belong to  $\mu_d=\set{\rev{e^{2\pi i k/d}}|\; 1\leq k\leq d}$. To study this spectrum we follow the $\Gamma$-space approach of \cite{gritschacher2019commuting}. This description allows us to prove that $k\mu_d$ is stably equivalent to $ku\wedge B\mu_d$. There is also a real version $ko_\sym$ constructed from commuting symmetric orthogonal matrices. We describe   low dimensional homotopy groups of these spectra.

\subsection{Classifying spaces}\label{sec:Classifying spaces}
Let $G$ be a topological group \revtwo{(locally compact and Hausdorff with a non-degenerate base point $1_G\in G$)}. An element $g\in G$ is said to be {\it $d$-torsion} if $g^d$ is the identity element $1_G$. We are interested in a space constructed from pairwise commuting $d$-torsion group elements.

\Def{\label{def:B-Z/d}
We define $B(\ZZ/d,G)$ to be the geometric realization of the simplicial space
$$
[n] \mapsto \Hom((\ZZ/d)^n,G)
$$
where $\Hom((\ZZ/d)^n,G)$ \rev{is the subspace of $G^n$ consisting of pairwise commuting $n$-tuples $(g_1,g_2,\cdots,g_n)$  such that $g_i^d = 1_G$ for all $1\leq i\leq n$.}  The simplicial structure is given by
$$
d_i(g_1,g_2,\cdots,g_n) = \left\lbrace 
\begin{array}{cc}
(g_2,\cdots,g_n) & i=0 \\
(g_1,\cdots, g_i g_{i+1},\cdots,g_n) & 0 <i < n \\
(g_1,g_2,\cdots,g_{n-1}) & i=n,
\end{array}
\right.
$$
and $s_j(g_1,g_2,\cdots,g_n)=(g_1,\cdots,g_j,1_G,g_{j+1},\cdots,g_n)$ for $0\leq j\leq n$.
}

\rev{
\Rem{{\rm
In general, for any cosimplicial group $\pi^\dt$ there is a classifying space $B(\pi,G)$ obtained by a similar construction, see \cite{okay2020commutative}. When $\pi^\dt$ is the level-wise free cosimplicial group $F^\dt$ then this construction gives the usual classifying space $BG$. If the level-wise abelianization $\ZZ^\dt$ is used then the resulting space is the {\it classifying space for commutativity} $B(\ZZ,G)$ \cite{AG15}. Mod-$d$ reduction in each level gives a cosimplicial group $(\ZZ/d)^\dt$ and we recover the construction given in Definition \ref{def:B-Z/d}. 
}}
}

\subsection{Stabilization}
Let $\CC^m$ denote the complex vector space of dimension $m$ with a canonical basis $\set{e_1,e_2,\cdots,e_m}$. Inclusion of the canonical basis vectors induces a map $\CC^m \to \CC^{m+1}$ and the union (colimit) along these inclusions is denoted by $\CC^{\infty}$. 
Let $U(m)$ denote the unitary group of $m\times m$ matrices. The stable unitary group $U$ is the union along the inclusions 
\begin{equation}\label{eq:stab}
U(m)\to U(m+1),\;\;\; A\mapsto \twobytwo{A}{0}{0}{1}.
\end{equation}
We write $B(\ZZ/d,U)$ for the union of $B(\ZZ/d,U(m))$ along the induced stabilization maps.

\rev{
\subsection{Spectra and $\Gamma$-spaces}  For basic properties of spectra we refer to \cite{adams1974stable,switzer2017algebraic}. A more recent exposition with applications to topological field theories can be found in \cite{beaudry2018guide}. 
A {\it spectrum} is a sequence $\set{E_n}_{n\geq 0}$ of pointed topological spaces together with pointed maps $\sigma_n: \Sigma E_n\to E_{n+1}$.
A morphism $f:E\to F$ between spectra consists of a sequence of pointed maps $\set{f_n:E_n\to F_n}_{n\geq 0}$ that commute with the structure maps $\sigma_n$.
Given a pointed topological space $X$ one can construct the suspension spectrum   $\Sigma^\infty X$ consisting of the $n$-fold suspensions $\set{\Sigma^n X}_{n\geq 0}$ where the structure maps $\sigma_n$ are given by the identity maps. The spectrum $\Sigma^\infty(S^0)$ is called the sphere spectrum and is denoted by $\sS$. 
There is a notion of homotopy for maps between spectra and one can talk about the set $[E,F]$ of homotopy classes of maps. Let $\Sigma^r E$ denote the shifted spectrum defined by $(\Sigma^r E)_n = E_{r+n}$ where   $r\in \ZZ$ (by convention $E_{r+n}$ is a point if $r+n<0$).
Homotopy groups of spectra are defined by
$$
\pi_r(E) = [\Sigma^r \sS, E].
$$
A cofiber sequence of spectra gives rise to a long exact sequence of homotopy groups. Given two spectra one can construct the smash product $E\wedge F$. When $F=\Sigma^\infty X$ we will write $E\wedge X$ for the corresponding smash product.
Spectra are used to define cohomology and homology theories. The $E$-cohomology and $E$-homology of $X$ are defined by
$$
E^r(X) = [\sS\wedge X,\Sigma^r E ]\;\;\;\text{ and }\;\;\; E_r(X) =  [\Sigma^r\sS, E\wedge X].
$$

Given a spectrum $E$ one can define a space  $\Omega^\infty E$  by taking the direct limit of the sequence of maps
$$
E_0 \xrightarrow{\omega_0} \Omega E_1  \xrightarrow{\Omega\omega_1} \cdots \xrightarrow{\Omega^{n-1}\omega_{n-1}} \Omega^nE_n \xrightarrow{\Omega^n\omega_n} \cdots
$$
where $\omega_n: E_n \to \Omega E_{n+1}$ is the adjoint of $\sigma_n$. The space $\Omega^\infty E$  has the structure of an {\it infinite loop space} \cite[\S 1.7]{adams1978infinite}. An infinite loop space can be delooped indefinitely, that is 
 one can define the spaces $\Omega^{-r} (\Omega^\infty E)$. 
Let $Q X$ denote the space $\Omega^\infty(\Sigma^\infty X)$. 
The $r$-th stable homotopy group of $X$ is defined by
$$
\pi_r^s(X) = \pi_r(QX) 
$$
which is also isomorphic to $\pi_r(\sS\wedge X)$.
 The $E$-cohomology of $X$ can also be defined as
$$
E^r(X) = [X, \Omega^{-r} (\Omega^\infty E)].
$$

}

\rev{
We are interested in spectra that come from  $\Gamma$-spaces. These objects are first introduced in \cite{segal1974categories}.
In this section we will mostly follow the 
 exposition given in \cite[\S 2]{gritschacher2019commuting}. For more details on $\Gamma$-spaces see \cite[\S 4]{bousfield1978homotopy} and \cite[App. B]{schwede2018global}.
} 
 Let $\catFin_*$ denote the category whose objects are pointed finite sets $k_+=\set{1,2,\cdots,k}\sqcup \set{+}$, $k\geq 0$, and morphisms are pointed set maps $\alpha:k_+ \to l_+$. Let $\catTop_*$ denote the category of pointed topological spaces. A {\it $\Gamma$-space} is a functor $F: \catFin_* \to \catTop_*$.
This can be extended to a functor $F:\catTop_*\to \catTop_*$ by a coend construction
\begin{equation}\label{eq:coend}
F(X) = \int^{k_+} F(k_+) \times X^k.
\end{equation}
\rev{More explicitly, the coend is the quotient of the disjoint union of $F(k_+)\times X^k$ over $k\geq 0$  under the equivalence relation generated by 
$$(F(\alpha)(z);x_1,\cdots,x_l)\sim (z;x_{\alpha(1)},\cdots ,x_{\alpha(k)})$$
where $\alpha$ runs over all pointed set maps $k_+\to l_+$ and $(z;x_1,\cdots,x_l)\in F(k_+)\times X^l$.
}
There is an assembly map $F(X)\wedge Y \to F(X\wedge Y)$ \rev{induced by $(z;x_1,\cdots,x_k)\wedge y\mapsto (z;x_1\wedge y,\cdots,x_k\wedge y)$.}
\rev{The assembly map can be used to associate a spectrum to a $\Gamma$-space. 
The spectrum associated to   $F$ is}  denoted by $F(\sS)$ \rev{and consists} of the spaces $\set{F(S^n)|\;n\geq 0}$ whose structure maps are induced by the assembly map \rev{$F(S^n)\wedge S^1 \to F(S^n\wedge S^1)$.} 

\begin{ex}\label{ex:gamma-spaces}  {\rm  \rev{We will encounter the following examples.}
\begin{enumerate}
\item   Let $\sS:\catFin_\ast \to \catTop_\ast$ denote the inclusion functor. This means that we regard $k_+$ as a pointed topological space with discrete topology. \rev{We can think of this functor as $\Hom_{\catFin_*}(1_+,-)$.}
The associated spectrum is the sphere spectrum and is simply denoted by $\sS$ \rev{\cite[\S 1]{schwede1999stable}}. 
  
\item \rev{Let $\CC^d$ denote the complex vector space with the canonical basis $e_1,e_2,\cdots,e_d$. The $\Gamma$-space $ku:\catFin_*\to \catTop_*$ is defined by}
$$
ku(k_+) = \coprod_{d_1,\cdots,d_k\in \NN}  \frac{L(\CC^{d_1}\oplus\cdots\oplus \CC^{d_k},\CC^\infty)}{U(d_1)\times \cdots \times U(d_k)}
$$
where $L(-,-)$ denotes the space of complex linear isometric embeddings between two complex inner product spaces. A point in this space is specified by a tuple $(V_1,\cdots,V_k)$ of pairwise orthogonal subspaces \rev{of $\CC^\infty$}. Given $\alpha:k_+\to l_+$ the map $ku(\alpha)$ is defined by
$$
(V_1,\cdots,V_k) \mapsto (\oplus_{i\in \alpha^{-1}(1)}V_i,\cdots, \oplus_{i\in \alpha^{-1}(l)}V_i ).
$$
The spectrum $ku(\sS)$ we obtain is the connective complex $K$-theory spectrum, which will be denoted simply by $ku$. \rev{There is a canonical morphism $\sS\to ku$ of $\Gamma$-spaces induced by the map $\sS(1_+) \to ku(1_+)=\coprod_{q\geq 0} \Gr_q(\CC^\infty)$ sending $1$ to the subspace $\Span{e_1} \subset \CC^\infty$ and the base point $+$ to the point $\Gr_0(\CC^\infty)$. }
There is a real version of this construction defined analogously but using $\RR$-vector spaces. The resulting spectrum is the connective real $K$-theory spectrum $ko$. \rev{ See \cite{gritschacher2019commuting} for an equivariant approach. }

\item Let $M$ be a commutative discrete monoid. \rev{There is an associated $\Gamma$-space denoted by $HM:\catFin_*\to \catTop_*$ where $HM(k_+)=M^k$} and for $\alpha:k_+\to l_+$ the map $\rev{HM}(\alpha)$ is defined by sending $(x_1,\cdots,x_k)$ to $(\sum_{j\in \alpha^{-1}(1)} x_j, \cdots,\sum_{j\in \alpha^{-1}(l)} x_j  )$. \rev{A}pplying $\Omega^\infty$  to the resulting spectrum $\rev{HM}(\sS)$ 
\rev{gives a space  homotopy equivalent to $\Omega BM$, also known as the group completion of $M$. This can be seen from Segal's delooping construction for $\Gamma$-spaces \cite[Pro. 1.4]{segal1974categories}.} 
In particular, we can consider the monoid $\NN$ and the associated spectrum $\rev{H}\NN(\sS)$. Since $\Omega^\infty \rev{H}\NN(\sS) \simeq \ZZ$ we see that this spectrum is equivalent to the Eilenberg-Maclane spectrum $H\ZZ$.
There is a map of $\Gamma$-spaces $\dim:ku\to \rev{H}\NN$ obtained by sending $(V_1,\cdots,V_k)$ to $(\dim(V_1),\cdots,\dim(V_k))$. \rev{A similar morphism of $\Gamma$-spaces $\dim:ku\to \rev{H}\ZZ$ can be obtained by replacing $\NN$ with the abelian group of integers. The canonical morphisms from $\sS$ fit into a commutative diagram of $\Gamma$-spaces}
$$
\begin{tikzcd}
\sS \arrow{r} \arrow{dr} & ku \arrow{d}{\dim} \\
& H\ZZ.
\end{tikzcd}
$$\rev{
The \revtwo{diagonal} arrow is induced by the inclusion map $1_+\to \ZZ$ defined by $1\mapsto 1$ and $+\mapsto 0$. }
\end{enumerate}
}
\end{ex}

\subsection{The spectrum}
Let $\mu_d\subset U(1)$ denote the subgroup generated by $e^{2\pi i/d}$.

\Pro{\label{pro:B-ku} Sending $(V_1,\cdots,V_k;\lambda_1,\cdots,\lambda_k)$, where $V_i$ are pairwise orthogonal finite-dimensional subspaces of $\CC^\infty$ and \rev{$\lambda_i\in (\mu_d)^n$}, to the $n$-tuple $(A_1,\cdots,A_n)$ of pairwise commuting unitary matrices, where \rev{$A_j$ acts on $V_i$} by multiplication with $\lambda_i^{(j)}$ and trivially on the complement of $V_1\oplus\cdots \oplus V_k$, induces a homeomorphism
$$
ku((\mu_d)^n) \xrightarrow{\cong} \Hom((\ZZ/d)^n,U).
$$
Moreover, this homeomorphism is compatible with the simplicial structures and induces a homeomorphism
$$
ku(B\mu_d)\xrightarrow{\cong} B(\ZZ/d,U).
$$
}
\Proof{The statements are proved in \cite{gritschacher2019commuting} when $\lambda_i^{(j)}\in U(1)^n$. These arguments still go through when $U(1)$ is replaced by the subgroup $\mu_d$.  
}

It is instructive to describe the inverse of the first homeomorphism. Let $(A_1,A_2,\cdots, A_n)$ be a tuple of pairwise commuting matrices in $U$ such that $(A_j)^d=I$ for $1\leq j\leq n$. These matrices are contained in $U(m)$ for some large enough $m$. 
We can simultaneously diagonalize these matrices
$$
\begin{pmatrix}
\lambda_1^{(1)} I_{d_1} &   &   &   \\
  & \lambda_2^{(1)} I_{d_2} &   &   \\
   &   & \ddots &   \\
  &   &   & \lambda_k^{(1)} I_{d_k} 
\end{pmatrix},\cdots,
\begin{pmatrix}
\lambda_1^{(n)} I_{d_1} &   &   &   \\
  & \lambda_2^{(n)} I_{d_2} &   &   \\
   &   & \ddots &   \\
  &   &   & \lambda_k^{(n)} I_{d_k} 
\end{pmatrix}
$$
such that $(\lambda_i^{(1)},\lambda_i^{(2)},\cdots,\lambda_i^{(n)})$ is distinct from $(\lambda_j^{(1)},\lambda_j^{(2)},\cdots,\lambda_j^{(n)})$ whenever $i\not=j$. Therefore $(A_1,A_2,\cdots,A_n)$ amounts to specifying a tuple $(V_1,V_2,\cdots,V_k)$ of  pairwise orthogonal finite dimensional subspaces $V_i\subset \CC^\infty$, $1\leq i\leq k$, together with the eigenvalues  $\lambda_i^{(j)}\in \mu_d$. Then the inverse map sends $(A_1,\cdots,A_n)$ to the class of $(V_1,\cdots,V_k;\lambda_1,\cdots,\lambda_k)$ in the coend construction \ref{eq:coend}.

Given a   pointed space $X$ and a $\Gamma$-space $F$ we write $F_X$ for the $\Gamma$-space defined by $F_X(k_+) = F(k_+ \wedge X)$. For $\alpha:k_+\to l_+$ the map $F_X(\alpha)$ is obtained by naturality of the coend construction.  A $\Gamma$-space $F$ is called {\it special} if the map 
$
F((k+l)_+) \to F(k_+) \times F(l_+)
$
induced by the projections $(k+l)_+ \to k_+$ and $(k+l)_+ \to l_+$ is a weak equivalence    for all $k_+,l_+$. A special $\Gamma$-space is called {\it very special} if $\pi_0 F(1_+)$ is an abelian group. \rev{(In general, $\pi_0 F(1_+)$ is a monoid since $F(1_+)$ is an $H$-space \cite[\S 1]{segal1974categories}.)}

\Lem{\label{lem:gamma} Let $X$ be a pointed space.
\begin{enumerate}
\item If $F$ is special   then $F_X$ is also special.

\item The natural map $F(\sS)\wedge X \to F_X(\sS)$ is a stable equivalence.

\end{enumerate}
}
\Proof{  Part (1) is implicitly mentioned in \rev{\cite[Proof of Theorem 4.4.]{bousfield1978homotopy}} and part (2) is  proved therein as Lemma 4.1. For a more recent exposition of the equivariant version of this statement  see \cite{schwede2018global} when $X$ has finitely many cells and  \cite{gritschacher2019commuting} for the general case.  }

\Def{{\rm
The spectrum $ku_{B \mu_d}(\sS)$ will be called the {\it commutative $d$-torsion $K$-theory spectrum} and will be denoted by  $k\mu_d$. The associated generalized cohomology theory will be referred to as the {\it commutative $d$-torsion $K$-theory}. \rev{We write $k\mu_d^n(X)$ to denote the $n$-th $k\mu_d$-cohomology of $X$ and  $k\mu_d(X)=k\mu_d^0(X)$ for simplicity of notation (not to be confused with the $\Gamma$-space evaluated at $X$)}
}}

\Pro{\label{pro:InfLoop}
The spectrum $k\mu_d$ is stably equivalent to $ku\wedge B\mu_d$ and the space  $\Omega^\infty k\mu_d$ is weakly equivalent to $B(\ZZ/d,U)$.
}
\Proof{ We modify the argument in \cite{gritschacher2019commuting} given for $B(\ZZ,U)$. Applying part (1) of the lemma to $F=ku$ and $X=B\mu_d$, and using the well-known fact that $ku$ is special we obtain that $ku_{B\mu_d}$ is special.  Moreover, $ku_{B\mu_d}$ is very special since
\begin{equation}\label{eq:B-ku-Gamma}
ku_{B\mu_d}(1_+) = ku(1_+\wedge B\mu_d) = B(\ZZ/d,U)
\end{equation}
and thus $\pi_0(ku_{B\mu_d}(1_+))=\pi_0 B(\ZZ/d,U)=0$. It is a general fact that if $F$ is very special then $\Omega^\infty F(\sS)\simeq F(1_+) $ \cite{segal1974categories}. Therefore $\Omega^\infty k\mu_d=\Omega^\infty ku_{B\mu_d}(\sS)\simeq ku_{B\mu_d}(1_+) \cong B(\ZZ/d,U)$. The equivalence $k\mu_d \simeq ku\wedge B\mu_d$ follows from part (2) of  Lemma \ref{lem:gamma}. 
}

\rev{As a consequence of this result we have
$$
k\mu_d^r(X) = [X, \Omega^{-r} B(\ZZ/d,U)].
$$
}

\Rem{\label{rem:diff-pi0}
{\rm
There is one important difference between $ku((\mu_d)^n)$ and $ku(U(1)^n)$ worth pointing out.  The former   is not an infinite loop space whereas the latter is since $U(1)^n$ is path connected. Note that $\pi_0 ku((\mu_d)^n)$ can be identified with $\Rep((\ZZ/d)^n,U)$, the union of the quotient spaces $\Hom((\ZZ/d)^n,U(m))/U(m)$ under the conjugation action of $U(m)$. 
 
Moreover, $\Rep((\ZZ/d)^n,U)\cong \rev{H}\NN((\mu_d)^n)$ and
the quotient map
$$
\Hom((\ZZ/d)^n,U)\to \Rep((\ZZ/d)^n,U)
$$ 
can be described using the map of $\Gamma$-spaces $\dim:ku\to \NN$; see \cite{gritschacher2019commuting}.  
For example, when $d=2$ we have that $\rev{H}\NN((\mu_2)^n)= \NN\wedge (\mu_2)^n$ where $\NN$ has $0$ as its base point and $(\mu_2)^n$ is based at the identity element. The set of path components is not an abelian group.
}
}

\subsection{Low dimensional homotopy groups}\label{sec:lowdim}
As a consequence of Proposition \ref{pro:InfLoop} \rev{the homotopy groups of $k\mu_d$ coincide} with $ku$-homology of $B\mu_d$. The groups $ku_*(B\mu_d)$ are computed in \cite[\S 3.4]{bruner2003connective}; see also  \cite{hashimoto1983connective}. In low degrees we have
\begin{equation}\label{eq:pi-Kcom}
\pi_r B(\ZZ/d,U) \cong\pi_r(k\mu_d) \cong\pi_r(ku \wedge B\mu_d) = \left\lbrace
\begin{array}{cc}
0 & r=0 \\
\ZZ/d & r=1 \\
0 & r=2.
\end{array}
\right.
\end{equation}

There is a commutative diagram
\begin{equation}\label{eq:Bmud}
\begin{tikzcd}
B\mu_d \arrow[r,hook] \arrow[rd,equal] & B(\ZZ/d,U) \arrow{d}{\det} \\
&B\mu_d
\end{tikzcd}
\end{equation}
which splits \rev{off} the $\ZZ/d$ in $\pi_1 B(\ZZ/d,U)$. 
The determinant map factors through the geometric realization of the simplicial set of connected components, denoted by $|\pi_0\Hom((\ZZ/d)^\dt,U)|$. Proposition \ref{pro:B-ku} implies that the connected components of $\Hom((\ZZ/d)^n,U)$ can be described as $\pi_0 ku((\mu_d)^n) = \rev{H}\NN((\mu_d)^n)$; see also Remark \ref{rem:diff-pi0}. 
Therefore we have
$$
\pi_0\Hom((\ZZ/d)^\dt,U) = \rev{H}\NN((\mu_d)^\dt)
$$
and the natural map $B(\ZZ/d,U)\to |\pi_0\Hom((\ZZ/d)^\dt,U)|$ is given by the geometric realization of
$$
ku((\mu_d)^\dt) \to \rev{H}\NN((\mu_d)^\dt)
$$
induced by the $\Gamma$-space map $\dim:ku\to \NN$ which sends a tuple of pairwise orthogonal subspaces $(V_1,V_2,\cdots,V_k)$ to their dimensions $(d_1,d_2,\cdots,d_k)$. Since $\NN$ is a special $\Gamma$-space we can apply Lemma \ref{lem:gamma} to obtain an equivalence
\begin{equation}\label{eq:N-mu-d-bullet}
|\rev{H}\NN((\mu_d)^\dt)|  \simeq \Omega^\infty (\rev{H}\NN(\sS) \wedge B\mu_d).
\end{equation}
Using the equivalence $\rev{H}\NN(\sS)\simeq H\ZZ$ we obtain the following.

\begin{pro}\label{pro:detfactor}
The determinant map factors as
$$
\begin{tikzcd}
B(\ZZ/d,U) \arrow{r} \arrow{rd}{\det} &{|\pi_0\Hom((\ZZ/d)^\dt,U)|} \arrow{d}\\
  & B\mu_d
\end{tikzcd}
$$
where the homotopy groups of the simplicial set of connected components is given by
$$
\pi_r |\pi_0\Hom((\ZZ/d)^\dt,U)| \cong \tilde H_r(B\mu_d,\ZZ).
$$
\end{pro}


The determinant map induces a  homomorphism
$$
\Det_*:k\mu_d(X) \to H^1(X,\ZZ/d)
$$
which splits as a consequence of the diagram in \ref{eq:Bmud}.
In general, since the homotopy groups of $k\mu_d$ are known we can compute $k\mu_d$-cohomology using the Atiyah--Hirzebruch spectral sequence \cite{adams1974stable}. The $E_2$-page of the spectral sequence is given by
\begin{equation}\label{eq:E2}
\rev{\tilde H}^p(X,\pi_{-q} k\mu_d) \Rightarrow k\mu_d^{\revtwo{p+q}}(X).
\end{equation}
One special case, for which the computation is easy, is when $X$ is a \rev{connected} $2$-dimensional CW complex. In this case the spectral sequence collapses in the $E_2$-page and $\Det_*$ becomes an isomorphism
\begin{equation}\label{eq:2dim}
k\mu_d(X) \cong H^1(X,\ZZ/d).
\end{equation}

\section{$C(d,m)$-cohomology}\label{sec:Cdm}

For each $m\geq 1$ we 
\revthr{construct} a spectrum, denoted by $C(d,m)$, obtained from the commutative $d$-torsion $K$-theory spectrum $k\mu_d$. 
In this section we compute the homotopy groups of $C(d,m)$ and describe the $C(d,m)$-cohomology of a space.  In \S \ref{sec:LCS} we will see that $C(d,m)$-cohomology informs us about operator solutions of linear constraint systems. These operator solutions play a significant role in quantum information theory.

\subsection{A quotient space}\label{sec:QuoSpa} Throughout this section \revtwo{we assume that the topological group $G$}
contains a central subgroup isomorphic to $\mu_d$. When $G=U(m)$ this will be the subgroup of $m\times m$ \rev{scalar matrices} with entries in $\mu_d$.

\Def{\label{def:bar-B-Z/d} Let $\bar B(\ZZ/d,G)$ denote the geometric realization of the simplicial space
$$
[n] \mapsto \Hom((\ZZ/d)^n,G)/\sim
$$
where the quotient relation identifies $(A_1,\cdots,A_n)$ with $(\alpha_1 A_1,\cdots,\alpha_n A_n)$ where $\alpha_i\in\mu_d$. Simplicial structure maps are similar to the ones given in Definition \ref{def:B-Z/d}.
} 

\rev{Let $\bar G$ denote the quotient group $G/\mu_d$.  The quotient space $\bar B(\ZZ/d,G)$ is a subspace of the classifying space $B\bar G$. Furthermore, there is a pull-back diagram}
$$
\begin{tikzcd}
B(\ZZ/d,G) \arrow{r}\arrow{d} & B G \arrow{d}\\
\bar B(\revtwo{\ZZ/d},G) \arrow{r} & B\bar G
\end{tikzcd}
$$
\rev{where the right-hand map is a fibration sequence with fiber $B\mu_d$. Therefore we obtain} a fibration sequence
$$
B\mu_d \xrightarrow{\Delta_G} B(\ZZ/d,G) \to \bar B(\ZZ/d,G)
$$
where the fiber inclusion is induced by $\mu_d\subset G$. By the classification of principal bundles this fibration is determined by a cohomology class $\gamma_G$ in $H^2(\bar B(\ZZ/d,G),\ZZ/d)$. 
When $G$ is the unitary group $U(m)$ we simply write $\Delta_m$ for the fiber inclusion and $\gamma_m$ for the cohomology class.   
The stabilization maps in  \ref{eq:stab} do not descend to $\bar B(\ZZ/d,U(m))$. However, we will construct a space which serves as a stabilization using methods from stable homotopy theory.  

\subsection{$C(d,m)$ spectrum}
We begin with a spectrum level description of  $\Delta_m$. 
For $m\geq 1$ let us introduce a map of $\Gamma$-spaces
\begin{equation}\label{eq:gamma-map}
\delta_m: \sS \to ku,
\end{equation}
induced by the map 
$$
1_+ \to \coprod_{\rev{q}\geq 0} \Gr_{\rev{q}}(\CC^\infty)
$$  
that sends  the element $1$ to the subspace $\CC^m=\Span{e_1,e_2,\cdots,e_m}$ and the base point $+$ to $\Gr_0(\CC^\infty)$.  This assignment determines all the other maps $\sS(k_+)\to ku(k_+)$ \rev{since  $\sS(-)\cong \Hom_{\catFin_*}(1_+,-)$.}

Let $\delta_{d,m}:\sS_{B\mu_d} \to ku_{B\mu_d}$ denote the $\Gamma$-space map induced by $\delta_m$ using the functoriality of the construction $F\mapsto F_X$. 
 The associated spectra maps will still be denoted by $\delta_m$ and $\delta_{d,m}$, respectively. 

Consider the cofiber sequence
\begin{equation}\label{eq:cofiber-C}
\sS_{B\mu_d}\rev{(\sS)} \xrightarrow{\delta_{d,m}} k\mu_d \to C(\delta_{d,m}).
\end{equation} 
\rev{For applications we are mainly interested in computing $[X,C(\delta_{d,m})]$ for a $2$-dimensional CW complex $X$. For a spectrum $X$ let $p_n:X\to P_n X$ denote the $n$-th Postnikov section \cite{lima1960}.} 

\Def{{\rm \label{def:barB-stable}
We define  \rev{$C(d,m)=P_2 C(\delta_{d,m})$, i.e.}  the   spectrum obtained from $C(\delta_{d,m})$ by killing the homotopy groups  of degree greater than $2$. We write $\bar B(d,m)$  for the associated infinite loop space $\Omega^\infty C(d,m)$.
}}   

To compute the homotopy groups of $C(d,m)$ we can use the cofiber sequence 
\begin{equation}\label{eq:cofib}
 \sS \wedge B\mu_d \xrightarrow{\delta_m\wedge \idy} ku \wedge B\mu_d \to C(\delta_{d,m})
\end{equation}
instead of \ref{eq:cofiber-C} since we have a commutative diagram of spectra 
\begin{equation}\label{eq:diag-Sp}
\begin{tikzcd}
\sS \wedge B\mu_d \arrow{r}{\sim} \arrow{d}{\delta_m\wedge \idy} & \sS_{B\mu_d}(\sS)\arrow{d}{\delta_{d,m}}\\
ku\wedge B\mu_d \arrow{r}{\sim} & ku_{B\mu_d}(\sS)
\end{tikzcd}
\end{equation}
as a consequence of part (2) of Lemma \ref{lem:gamma}.

\Lem{\label{lem:Cdm-homotopy}The homotopy groups of $C(d,m)$ fit into an exact sequence \rev{of abelian groups}
\begin{equation}\label{eq:homotopy-C}
0 \to \pi_2 C(d,m) \to \ZZ/d \xrightarrow{\phi}  \ZZ/d \to  \pi_1 C(d,m)\to 0
\end{equation}
where $\phi$ is given by multiplication with $m$.}
\begin{proof}
\rev{The exact sequence in \ref{eq:homotopy-C} is obtained from the homotopy exact sequence associated to \ref{eq:cofib} and using the homotopy groups of $ku\wedge B\mu_d$ given in \ref{eq:pi-Kcom} together with the isomorphism $\pi_1^s(B\mu_d)\cong \ZZ/d$. It remains to show that $\phi:\ZZ/d\to \ZZ/d$ is given by multiplication with $m$.
The $\Gamma$-space map $\revtwo{\delta_{m}}:\sS \to ku$ fits into a commutative diagram of $\Gamma$-spaces}
$$
\begin{tikzcd}
\sS \arrow[d,"\revtwo{\delta_{m}}"'] \arrow{r} & H\ZZ \arrow{d}{\times m} \\
ku \arrow{r}{\dim} & H\ZZ
\end{tikzcd}
$$
\rev{where the top map is the canonical map of $\Gamma$-spaces (Example \ref{ex:gamma-spaces} part (3)). The right-hand map is induced by  $\ZZ\xrightarrow{\times m}\ZZ$, the multiplication with $m$ map.   Replacing each $\Gamma$-space $F$ in the diagram with $F_{B\mu_d}$, looking at the associated spectrum and using part (2) of Lemma \ref{lem:gamma} we obtain a diagram of spectra}
\begin{equation}\label{eq:diag-spectra}
\begin{tikzcd}
\sS \wedge B\mu_d \arrow[d,"\revtwo{\delta_{m}}\wedge \idy"'] \arrow{r} & H\ZZ \wedge B\mu_d \arrow{d}{(\times m)\wedge \idy} \\
ku \wedge B\mu_d \arrow{r} & H\ZZ \wedge B\mu_d
\end{tikzcd}
\end{equation}
\rev{Applying $\pi_1$ to this diagram gives a commutative diagram of abelian groups}
$$
\begin{tikzcd}
\pi_1^s(B\mu_d) \arrow{d}{\phi} \arrow{r}{\cong} & H_1(B\mu_d,\ZZ) \arrow{d}{\times m} \\
ku_1(B\mu_d) \arrow{r}{\cong} & H_1(B\mu_d,\ZZ)
\end{tikzcd}
$$
\rev{Therefore after identifying these groups with $\ZZ/d$ we see that $\phi$ is given by   multiplication with $m$.
}

\end{proof}

\rev{
We can extend the diagram \ref{eq:diag-spectra} by composing the horizontal maps with the first Postnikov section}
\begin{equation}\label{eq:diag-spectra-extended}
\begin{tikzcd}
\sS \wedge B\mu_d \arrow[d,"\revtwo{\delta_{m}}\wedge \idy"'] \arrow{r} & H\ZZ \wedge B\mu_d \arrow{d}{(\times m)\wedge \idy} \arrow{r}{p_1} & \Sigma H\ZZ/d  \arrow{d}{\times m} \\
ku \wedge B\mu_d \arrow{r} & H\ZZ \wedge B\mu_d \arrow{r}{p_1} &\Sigma  H\ZZ/d
\end{tikzcd}
\end{equation}
\rev{

\Cor{\label{cor:Cdm-cofiber}There is a cofiber sequence
\begin{equation}\label{eq:cofib-Cdm}
\Sigma H\ZZ/d \xrightarrow{\times m} \Sigma H\ZZ/d \to C(d,m).
\end{equation}
\revtwo{Moreover, we have an equivalence $ \Sigma H\pi_1(C(d,m))\vee \Sigma^2 H\pi_2(C(d,m))\to C(d,m)$.}
}
\Proof{
\revtwo{The cofiber sequence is a consequence of}
Lemma \ref{lem:Cdm-homotopy}. 
\revtwo{This sequence  presents $C(d,m)$ as a $H\ZZ$-module spectrum since $\times m$ is a morphism of $H\ZZ$-modules.
For a $H\ZZ$-module spectrum $E$ there is an equivalence $\vee_{n\in \ZZ} \Sigma^n H\pi_n(E)\to E$ of spectra \cite[Pro. 5.3]{casacuberta2005homotopical}, which gives the splitting.
}
}}

\subsection{$C(d,m)$-cohomology} 
Let us introduce notation for the abelian groups corresponding to the kernel and the cokernel of the exact sequence in \ref{eq:homotopy-C} 
$$
0\to (\ZZ/d)_m \xrightarrow{i_m} \ZZ/d \xrightarrow{\times m} \ZZ/d \xrightarrow{\pi_m} \frac{\ZZ/d}{m\ZZ/d} \to 0.
$$
For a group homomorphism $h:A\to B$ we write $h_*:H^n(X,A)\to H^n(X,B)$ for the change of coefficients map.   Note that both the kernel and the cokernel are isomorphic to $\ZZ/\gcd(d,m)$.

\rev{The $C(d,m)$-cohomology of a pointed space $X$ is defined by the pointed homotopy classes of maps
$$
C(d,m)^r(X) = [X, \Omega^{-r} \bar B(d,m)]
$$
and we simply write $C(d,m)(X)=[X,\bar B(d,m)]$ when $r=0$.
}

\begin{thm}\label{thm:Thm-stable} 
Let $X$ be a connected 
CW complex. There is a commutative diagram
\begin{equation}\label{eq:diag-thm}
\begin{tikzcd}
& H^2(X,(\ZZ/d)_m) \arrow[d,hook]\arrow[r,"(i_m)_*"] &  H^2(X,\ZZ/d) \arrow[d,equal] \\
k\mu_d(X) \arrow[d] \arrow{r}{\zeta} & C(d,m)(X) \arrow[d,two heads] \arrow{r}{\cl} & H^2(X,\ZZ/d) \\ 
H^1(X,\ZZ/d) \arrow[r,"(\pi_m)_*"] & H^1(X,\frac{\ZZ/d}{m\ZZ/d}) &
\end{tikzcd}
\end{equation}
where  $\cl(f)=f^*(\gamma_m^\sS)$, the image of $\zeta$ is contained in the kernel of $\cl$, and the middle column is an exact sequence 
\revtwo{which (non-canonically) splits as follows
$$
C(d,m)(X) \cong H^1(X,\frac{\ZZ/d}{m\ZZ/d}) \oplus H^2(X,(\ZZ/d)_m) .
$$}
\end{thm}
\begin{proof}
\rev{
The diagram of spectra given in \ref{eq:diag-spectra-extended} extends below to the cofibers. Shifting the resulting diagram of cofibrations gives the following diagram of spectra
}
$$
\begin{tikzcd}
k\mu_d \arrow{r} \arrow{d}{c} & \Sigma H\ZZ/d \arrow{d}\\
C(\delta_{d,m}) \arrow{r}{p_2}\arrow{d} & C(d,m) \arrow{d}{\gamma} \\
\Sigma (\sS \wedge B\mu_d) \arrow{r} & \Sigma^2 H\ZZ/d
\end{tikzcd}
$$ 
The exactness claim about the middle row \rev{of \ref{eq:diag-thm}}  is obtained by evaluating the sequence $k\mu_d \xrightarrow{p_2 c} C(d,m) \xrightarrow{\gamma} \Sigma^2 H\rev{\ZZ/d}$ at the space $X$. \rev{Applying \revtwo{to} $\gamma:C(d,m)\to\Sigma^2 H\ZZ/d$ the functor $\Omega^\infty$ gives a map  $\Omega^\infty\gamma:\bar B(d,m)\to B^2 \mu_d$ which represents the cohomology class $\gamma_m^\sS \in H^2(\bar B(d,m),\ZZ/d)$. Therefore for a map $f:\sS\wedge X\to C(d,m)$ the cohomology class $\cl(f)$, which is represented by  the composition $\gamma f$, coincides with $f^*(\gamma_m^\sS)$.
}
\rev{The horizontal maps $(\pi_m)_*$ and $(i_m)_*$ are obtained by comparing the Atiyah--Hirzebruch spectral sequences and the commutativity of the squares follow from the naturality of the spectral sequences.} The 
\revtwo{homomorphism} $k\mu_d(X) \to H^1(X,\ZZ/d)$ 
is  the edge homomorphism in the Atiyah--Hirzebruch spectral sequence for $k\mu_d$-cohomology of $X$. The splitting \revtwo{is given by the spectrum level splitting described in} Corollary \ref{cor:Cdm-cofiber}.
\end{proof}   

We end this section by considering the associated infinite loop space $\bar B(d,m)= \Omega^\infty C(d,m)$ and its relation to the unstable spaces $\bar B(\ZZ/d,U(m))$. \revtwo{We will identify  $\ZZ/d$ with $\mu_d$ via the isomorphism $1\mapsto \omega$.} Recall that the fibration
$$
B\revtwo{(\ZZ/d)} \xrightarrow{\Delta_m} B(\ZZ/d,U(m)) \to \bar B(\ZZ/d,U(m))
$$
is determined by a cohomology class $\gamma_m \in H^2(\bar B(\ZZ/d,U(m)),\ZZ/d )$. Applying \revtwo{the functor} $\Omega^\infty$  to the cofiber sequence in \ref{eq:cofib-Cdm} 
gives a (homotopy) fiber sequence
\begin{equation}\label{eq:barB-fibration}
B\revtwo{(\ZZ/d)} \xrightarrow{B(\times m)} B\revtwo{(\ZZ/d)} \to \bar B(d,m)
\end{equation}
which is determined by a cohomology class $\gamma_m^\sS \in H^2(\bar B(d,m),\ZZ/d)$. 
\revtwo{Observe that the determinant map $\det:B(\ZZ/d,U(m))\to B(\ZZ/d)$ descends to a map $\overline\det:\bar B(\ZZ/d,U(m)) \to B(\frac{\ZZ/d}{m\ZZ/d})$.} 


\begin{lem}\label{lem:stablecan}
\begin{enumerate}
\item \rev{There \revtwo{are maps}  of (homotopy) fibrations}
 \begin{equation}\label{eq:diag-Q-fib}
\begin{tikzcd}
B\revtwo{(\ZZ/d)} \arrow[r,equal] \arrow{d}{\Delta_m} & B\revtwo{(\ZZ/d)} \arrow{d}{B(\times m)} \arrow[r,"B(\times m)"] & B(m\ZZ/d) \arrow{d} \\
B(\ZZ/d,U(m)) \arrow[r,"\det"] \arrow{d}& B\revtwo{(\ZZ/d)} \arrow{d} \arrow[r,equal] & B(\ZZ/d)  \arrow[d,"B\pi_m"] \\
\bar B(\ZZ/d,U(m)) \arrow[r,"\bar\iota_m"] \arrow[rr,bend right,"\overline\det"] & \bar B(d,m) \arrow[r,"p_1"] & B\left(\frac{\ZZ/d}{m\ZZ/d}\right) .
\end{tikzcd}
\end{equation}
where the diagram commutes up to homotopy \revtwo{and  $p_1\bar\iota_m \simeq \overline\det$}.
\item \revtwo{The cohomology class} $\gamma_m^\sS$ maps to $\gamma_m$ under the induced map $(\bar \iota_m)^*$. \revtwo{Moreover, the maps representing these cohomology classes fit into a homotopy commutative diagram} 
\begin{equation}\label{eq:diag-bar-gamma_m}
\begin{tikzcd}
\bar B(\ZZ/d, U(m)) \arrow[dd,bend right=60,"\gamma_m"'] \arrow[r,"\bar\iota_m"] \arrow[d,"\bar\gamma_m"] & \bar B(d,m) \arrow[d,"{\bar \gamma_m^\sS}"] \arrow[dd,bend left=60,"\gamma_m^\sS"] \\
B^2(\ZZ/d)_m \arrow[r,equal] \arrow[d,"B^2i_m"] & B^2(\ZZ/d)_m \arrow[d,"B^2i_m"] \\
B^2(\ZZ/d) \arrow[r,equal] & B^2(\ZZ/d)
\end{tikzcd}
\end{equation}
\revtwo{where $B^2(-)=B(B(-))$.}
\end{enumerate}
\end{lem}
\begin{proof}
\revtwo{Part (1):} Let $G$ be a topological group and $X$ be a $G$-space. The homotopy quotient of $X$, also known as the Borel construction \cite[\S 2.2]{adem2002topics}, will be denoted by $X\hquo G$. A $G$-map $X\to Y$ between two $G$-spaces induces a map $X\hquo G \to Y\hquo G$ between the homotopy quotients. We will apply this idea to construct the map $\bar\iota_m$ in diagram \ref{eq:diag-Q-fib} which  will make the lower part of the diagram commute up to homotopy. Since $\revtwo{\ZZ/d}$ is an abelian group the classifying space $B\revtwo{(\ZZ/d)}=|\revtwo{(\ZZ/d)}^\dt|$ is a topological \revtwo{abelian} group. The group structure is induced by the levelwise \revtwo{addition} maps 
$\set{m_n:\revtwo{(\ZZ/d)}^n\times \revtwo{(\ZZ/d)}^n \to \revtwo{(\ZZ/d)}^n  }_{n\geq 0}$ defined by 
$$m_n((\lambda_1,\lambda_2,\cdots,\lambda_n ) ,(\lambda'_1,\lambda_2',\cdots,\lambda'_n)) = \revtwo{(\lambda_1+\lambda'_1,\lambda_2+\lambda'_2,\cdots,\lambda_n+\lambda_n')}.$$ 
More precisely $\revtwo{(\ZZ/d)}^\dt$ is a simplicial abelian group. Moreover, $B(\ZZ/d,U(m))$ and $B\revtwo{(\ZZ/d)}$ are $B\revtwo{(\ZZ/d)}$-spaces. Both spaces are geometric realizations of simplicial spaces and the action is given by the geometric realization of a levelwise (simplicial) action. For the first one the action on the space of $n$-simplices is  given by 
$$
(\lambda_1,\lambda_2,\cdots,\lambda_n)\cdot (A_1,A_2,\cdots,A_n) \mapsto (\revtwo{\omega^{\lambda_1}}A_1,\revtwo{\omega^{\lambda_2}}A_2,\cdots,\revtwo{\omega^{\lambda_n}}A_n)
$$
and for the second space the action is given by
$$
(\lambda_1,\lambda_2,\cdots,\lambda_n)\cdot (\lambda'_1,\lambda'_2,\cdots,\lambda'_n) \mapsto (\revtwo{m\lambda_1}+\lambda'_1,\revtwo{m\lambda_2}+\lambda'_2,\cdots,\revtwo{m\lambda_n}+\lambda'_n).
$$
With these actions $\det:B(\ZZ/d,U(m))\to B\revtwo{(\ZZ/d)}$ is $B\revtwo{(\ZZ/d)}$-equivariant, Moreover, the base spaces in diagram \ref{eq:diag-Q-fib} are homotopy equivalent to the corresponding homotopy quotients.
Then the map $\bar\iota_m$ is obtained as the composite
$$
\bar B(\ZZ/d,U(m))\simeq B(\ZZ/d,U(m)) \hquo B\revtwo{(\ZZ/d)} \to B\revtwo{(\ZZ/d)}\hquo B\revtwo{(\ZZ/d)}  \simeq \bar B(d,m).
$$
\revtwo{A similar argument applies to the right-hand square and produces the map $p_1$  between the homotopy quotients. The comparison maps between the homotopy quotients and the ordinary quotients give a commutative diagram}
$$
\begin{tikzcd}
B(\ZZ/d,U(m)) \hquo B(\ZZ/d) \arrow[d,"\sim"] \arrow[r,"p_1\bar\iota_m"] &  B(\ZZ/d) \hquo B(m\ZZ/d)  \arrow[d,"\sim"]\\ 
\bar B(\ZZ/d,U(m)) \arrow[r,"\overline\det"] &  B\left(\frac{\ZZ/d}{m\ZZ/d}\right)
\end{tikzcd}  
$$

\revtwo{Part (2): Delooping the left-hand square of diagram (\ref{eq:diag-Q-fib}) shows that the maps $\gamma_m:\bar B(\ZZ/d,U(m))\to B^2(\ZZ/d)$ and $\gamma_m^\sS:\bar B(d,m)\to B^2(\ZZ/d)$ representing the corresponding cohomology classes fit into the outer diagram in  (\ref{eq:diag-bar-gamma_m}). In particular, this implies that $\bar\iota_m^*(\gamma_m^\sS)=\gamma_m$. Since up to homotopy $\gamma_m^\sS$ is the homotopy fiber of the map $B^2(\ZZ/d) \xrightarrow{B^2\times m} B^2( \ZZ/d)$ it factors as $\bar\gamma_m^\sS: \bar B(d,m)\to B^2 (\ZZ/d)_m$. A similar factorization applies to $\gamma_m$ since it is obtained as the pull-back of $\gamma_m^\sS$ along $\bar \iota_m$ giving us the homotopy commutative diagram (\ref{eq:diag-bar-gamma_m}).   
}

\end{proof}

\begin{cor}\label{cor:iota-bar-alternative} \revtwo{
The following diagram commutes up to homotopy}
$$
\begin{tikzcd}
\bar B(\ZZ/d,U(m)) \arrow[r,"{\bar \iota_m}"] \arrow[dr,"{\overline\det\times \bar \gamma_m}"'] & \bar B(d,m) \arrow[d,"p_1\times \bar \gamma_m^\sS","\sim"'] \\
& B\left(\frac{\ZZ/d}{m\ZZ/d}\right)\times B^2(\ZZ/d)_m  
\end{tikzcd}
$$ 
\end{cor} 
\begin{proof}
\revtwo{
Both $(p_1\times \bar\gamma_m^\sS)\bar\iota_m$ and $\overline\det\times \bar \gamma_m$ correspond to the same class in 
$$H^1(\bar B(\ZZ/d,U(m)),\frac{\ZZ/d}{m\ZZ/d})\times H^2(\bar B(\ZZ/d,U(m)),(\ZZ/d)_m)$$
 since $p_1\bar\iota_m\simeq \overline\det$ and $\bar\gamma_m^\sS\bar\iota_m\simeq \bar\gamma_m$  by Lemma \ref{lem:stablecan}. Therefore $(p_1\times \bar\gamma_m^\sS)\bar\iota_m\simeq\overline\det\times \bar \gamma_m$. Note that $p_1\times \bar\gamma_m^\sS$ is a homotopy inverse of the homotopy equivalence obtained by applying $\Omega^\infty$ to the equivalence in Corollary \ref{cor:Cdm-cofiber}.
}
\end{proof}

 
  

\section{Operator solutions of linear constraint systems}\label{sec:LCS}

Linear constraint systems arise in quantum information theory in the context of non-local games. Such games are played among a referee and two players where each player aims to   win the game  by satisfying a fixed set of rules. For some games if the players use  quantum resources, such as  entangled quantum states and  quantum measurements, then they can increase their likelihood of winning the game. Other than their applications in quantum information theory, linear constraint systems have found applications in resolving problems in the theory of operator algebras such as Tsirelson problem \cite{slofstra2019tsirelson} and Connes embedding conjecture \cite{ji2020mip}. In this section we study operator solutions of linear constraint systems by using the generalized cohomology theory, $C(d,m)$-cohomology, introduced in \S \ref{sec:Cdm}. We show that operator solutions of linear constraint systems correspond to classes in $C(d,m)$-cohomology.
 The paradigmatic example of a linear constraint system constructed by Mermin \cite{Mermin}\rev{, see also \cite{mermin1990simple,peres1990incompatible},} gives rise to a 
  non-trivial class in the $C(2,2^n)$-cohomology of a torus for  $n\geq 2$.   
This connection to stable homotopy theory opens up a new direction in the study of linear constraint systems. In this respect stable homotopy theory plays a similar role as it does
in the classification of topological quantum phases \cite{kitaev2013classification}; see also \cite{marcolli2019gamma} for applications of stable homotopical methods to quantum information theory.

\subsection{Linear constraint systems}
A {\it linear constraint system}  is specified by a system of linear equations $Mx=b$ for some $r\times c$ matrix $M$ with entries in $\ZZ/d$ \rev{and $b\in (\ZZ/d)^r$}.  
We say that a linear constraint system has an {\it operator solution}   if there exists a collection of $m\times m$-unitary matrices $A_i$, $1\leq i\leq c$, such that
\begin{enumerate}
\item $(A_i)^d$ is the identity matrix $I_m$ for all $1\leq i\leq c$,
\item $A_{i} A_j=A_j A_{i}$ whenever $M_{ki}$ and $M_{kj}$ are both non-zero for some $1\leq k\leq r$,
\item $A_1^{M_{k1}} A_2^{M_{k2}}\cdots A_{c}^{M_{{kc}}} = \omega^{b_k}I_m$, where $\omega=e^{2\pi i /d}$, for all $1\leq k\leq r$.
\end{enumerate}
When $m=1$ we call such a solution a {\it scalar solution}. In the physics literature an operator solution is usually called a {\it quantum solution} and a scalar solution is called a {\it classical solution}. A linear constraint system which admits no classical solutions is called {\it contextual}; otherwise it is called {\it non-contextual}. Note that in this paper we restrict our attention to operator solutions over finite-dimensional Hilbert spaces. The finiteness restriction can be removed for a more general discussion of the subject. For basic properties of linear constraint systems we refer to \cite{cleve2014characterization,Slo,qassim2020classical,
okay2020homotopical}.

\subsection{Topological description}
A linear constraint system can be formulated using hypergraphs. The data of a linear constraint system can be turned into a pair $(\hH,\tau)$ where $\hH=(V,E,\epsilon)$ is a hypergraph with \rev{a finite} vertex set $V$, \rev{a finite} edge set $E$ and an incidence weight $\epsilon$; and $\tau$ is a function $E\to \ZZ/d$. 
More concretely, let $\hH$ denote the hypergraph with $V=\set{v_1,v_2,\cdots,v_c}$, $E=\set{e_1,e_2,\cdots ,e_r}$ where $e_k=\set{v_i|\; M_{ki}\neq 0}$, and   $\epsilon_{e_k}(v_i)=M_{ki}$. 
The function $\tau$ is defined by $\tau(e_k)=b_k$. An operator solution can be regarded as a function $T:V\to U(m)$ where $T(v_i)=A_i$.

\Rem{{\rm
As in \S \ref{sec:QuoSpa} let $G$ be a group which contains a central subgroup isomorphic to $\mu_d$. We can consider solutions over $G$ instead of $U(m)$. Such an operator solution will be denoted by a function $T:V\to G$ where the group elements $\set{T(v)|v\in V}$ satisfy the $d$-torsion (1), commutativity (2), and linear constraint (3) conditions listed above.
}}

We define a chain complex associated to the hypergraph
$$
C_*(\hH):C_2 \stackrel{\partial}{\to} C_1 \stackrel{0}{\to} C_0
$$
where 
$$C_0=\mathbb{Z}/d,\;\;C_1=\mathbb{Z}/d[V],\;\;C_2=\mathbb{Z}/d[E],\;\;\;\;\partial[e]=\sum_{v\in e} \epsilon_e(v) [v].$$
 There is a corresponding cochain complex $C^*(\hH)$. The function $\tau$ can be regarded as a \rev{$2$-cocycle}. We write $[\tau]$ for its cohomology class.

\rev{Let $X$ be a CW complex.  
The set of $n$-cells will be denoted by $X_n$ and $X^n$ will denote the $n$-skeleton. For each $n$-cell there is an attaching map $\varphi^n: \partial D^n \to X^{n-1}$ and a characteristic map $\Phi^n:D^n \to X$. The exact sequence of   homotopy groups associated to the pair $(X^2,X^1)$ of complexes is given by
\begin{equation}\label{eq:fundamental-seq}
0\to \pi_2(X^2) \to \pi_2(X^2,X^1) \xrightarrow{\partial} \pi_1(X^1) \to \pi_1(X^2) \to 0.
\end{equation}
If $X$ has a single $0$-cell then the fundamental group is   generated by the homotopy classes $[\Phi^1]$ of the characteristic maps of $1$-cells. The relative homotopy group $\pi_2(X,X^1)$ is generated\revtwo{, up to the 
action of $\pi_1(X^1)$,} by the homotopy classes \revtwo{of the characteristic maps $[\Phi^2]$} of $2$-cells \cite[Chapter II, \S 2.1]{TwoDim}. 
This observation applies to the quotient space   $\bar X=X/X_0$   obtained by identifying the $0$-cells in $X$. Let $q:X\to \bar X$ denote the quotient map. We will write $q^{n}:X^n \to \bar X^n$ to denote the induced map between the $n$-skeletons. 
The characteristic maps \revtwo{of $\bar X$} for $n\geq 1$ are given by the composites $\bar \Phi^n: D^n \xrightarrow{\Phi^n} X \xrightarrow{q} \bar X$.
}
 
\begin{defn} \label{def:topological-real} {\rm A \emph{topological realization} for the hypergraph $\hH$ is a connected $2$-dimensional CW complex \rev{$X=X(\hH)$ that satisfies the following properties:
\begin{enumerate}
\item There are isomorphisms of sets: $X_1\cong V$ and $X_2\cong E$. The attaching (characteristic) maps   corresponding to  $1$-cells and $2$-cells will be labeled by $V$ and $E$, respectively.

\item For each $e\in E$ the image of $\varphi^2_e: \partial D^2 \to X^1$ is contained in the union of the images of the characteristic maps $\Phi_v^1:D^1 \to X$ where $v\in e$.

\item  For each $e\in E$ the boundary map $\partial: \pi_2(\bar X^2,\bar X^1) \to \pi_1(\bar X^1)$  in \ref{eq:fundamental-seq} satisfies
$$
\partial[\revtwo{\bar \Phi}_e^2] = \prod_{v\in e} [\bar \Phi_v^1]^{\epsilon_e(v)}
$$ 
where the product is with respect to some ordering of the set $e$.  
\end{enumerate} }
}
\end{defn} 

\begin{rem}{\rm 
\rev{The notion of topological realization introduced above is more restricted than the one introduced in \cite[Def. 6.1]{okay2020homotopical}.
Let $C_*(X)$ denote the chain complex $X$ over $\ZZ/d$. Definition \ref{def:topological-real} induces a morphism}  of chain complexes $f_*:C_*(X) \to C_*(\hH)$ such that $f_1$ and $f_2$ are \rev{isomorphisms}, i.e.,
$$
\begin{tikzcd}
 C_2(X) \arrow{r}{\partial_2} \arrow[d,"f_2","\cong"'] &  C_1(X)  \arrow{r}{\partial_1}  \arrow[d,"f_1","\cong"'] & C_0(X)  \arrow[d,"f_0"] \\
 \ZZ/d[E] \arrow{r}{\partial}   &  \ZZ/d[V]  \arrow{r}{0}  & \ZZ/d   
\end{tikzcd}
$$

}
\end{rem}
 
\begin{construction}\label{con:fT}
{\rm  \rev{Let $(\hH,\tau)$ be a LCS and $T:V\to G$ be an operator solution over $G$. For each edge $e\in E$ in the hypergraph define the  subgroup $A(e)\subset G$ \revtwo{(which is regarded as a discrete group)} generated by $\mu_d \cup \set{\rev{T(v)}|\;v\in e}$ and let $\bar A(e) = A(e)/\mu_d$. Let $X=X(\hH)$ be a topological realization for $\hH$. 
W}e  construct a map
$$
f_T:X \to \bar B(\ZZ/d,G),
$$
which is defined up to homotopy, as follows
\begin{enumerate}
\item send each $0$-cell in $X_0$ to the unique \rev{$0$-cell} of $\bar B(\ZZ/d,G)$,
\item send the $1$-cell labeled by $v\in X_1$ to the $1$-cell labeled by $[T(v)]$, the equivalence class of $T(v)$ under multiplication with elements in $\mu_d$,

\item \rev{by parts (2) and (3) of Definition \ref{def:topological-real}} the boundary of a $2$-cell labeled by $e\in X_2$ maps to a contractible loop in the subspace $B\bar A(e) \subset \bar B(\ZZ/d,G)$; extend this map to the interior of the disk by choosing a contracting homotopy   that lies in $B\bar A(e)$.  
\end{enumerate}
}
\end{construction}

\begin{rem}
{\rm \revtwo{Part (3) of Construction \ref{con:fT} requires some explanation. Let $Y=\bar B(\ZZ/d,G)$.} Observe that $f_T$ factors through the quotient map $q:X\to \bar X$ since all the $0$-cells are identified. \revtwo{Let $\bar f:\bar X\to Y$ be such that $\bar f q=f_T$. This map induces a commutative diagram (see \ref{eq:fundamental-seq})}
$$
\begin{tikzcd}
\pi_2(\bar X, \bar X_2) \arrow[r,"\partial"] \arrow[d,"\bar f_2"] & \pi_1(\bar X^1) \arrow[d,"\bar f_1"] \\
\pi_2(Y,Y^1) \arrow[r,"\partial"] & \pi_1(Y^1)
\end{tikzcd}
$$
\revtwo{For a $2$-cell labeled by $e\in X_2$ we have 
$$\partial \bar f_2[\bar\Phi_e^2] = \bar f_1 \partial [\bar\Phi_e^2]= \prod_{v\in e} \bar f_1[\bar\Phi_v^1]^{\epsilon_e(v)} = \prod_{v\in e} [T(v)]^{\epsilon_e(v)}=1  $$
as a consequence of the relation satisfied by the group elements $\set{T(v)}_{v\in e}$} (as part of the definition of a LCS). \revtwo{Therefore} the homotopy class of \revtwo{the image of} the boundary of the $2$-cell \revtwo{$\bar\Phi_e^2$} is contractible. 
Moreover,
 any  two choices of a contracting homotopy extending the map on the boundary of a $2$-cell of $X$ are homotopic to each other since the image lies inside the subspace $B\bar A(e)$ \revtwo{(by part (2) of Definition \ref{def:topological-real})}, whose homotopy groups above degree $2$ vanishes.  Therefore the map $f_T$ is unique up to homotopy. 
}
\end{rem}

Let $[(X,x_0),(Y,y_0)]$ denote the set of pointed homotopy classes of maps between two based spaces. We will suppress the base points and simply write $[X,Y]$. This should not result in any confusion since in this paper we do not consider the set of unpointed homotopy classes of maps.

\Pro{\label{pro:LCS-basics} Let $(\hH,\tau)$ be a linear constraint system.
\begin{enumerate}
\item $(\hH,\tau)$ has a scalar solution if and only if $[\tau]=0$ in $H^2(C(\hH))$, and thus, in the second cohomology group of any topological realization.

\item If $T$ is an operator solution for $(\hH,\tau)$ then $f_T^*(\gamma_G)=[\tau]$ for any map $f_T$ constructed using the operator solution (Construction \ref{con:fT}).

\item If $(\hH,\tau)$ has an operator solution $T$ and a topological realization $X$ such that $f_T$ induces the trivial map between the fundamental groups   then $(\hH,\tau)$   has a scalar solution.

\end{enumerate}
}
\begin{proof}
\rev{Part} (1) follows from the definition of the chain complex, see also \cite{Coho}.  Part (2) follows from the observation \revtwo{that $\gamma_G$ is the image of the identity map in $H^1(\mu_d,\ZZ/d)\cong \Hom(\ZZ_d,\ZZ_d)$ under the transgression map of the fibration $B\mu_d \to B(\ZZ/d,G)\to \bar B(\ZZ/d,G)$. In effect the transgression map is computed using the connecting homomorphism $\delta:H^1(B\mu_d,\ZZ/d)\to H^2((B(\ZZ/d,G),B\mu_d),\ZZ/d)$ in the cohomology long exact sequence of the pair $(B(\ZZ/d,G),B\mu_d)$; see \cite[page 186]{mccleary2001user}. Let $\phi:\mu_d\to \ZZ/d$ be the $1$-cocycle representing the identity map in cohomology. By definition of the connecting homomorphism $\delta[\phi]$ will be represented by the coboundary of a lift of $\phi$ to a $1$-cochain on $B(\ZZ/d,G)$. To   compute the value of the pull-back $f_T^*(\gamma)$ of the resulting $2$-cocycle  $\gamma$ on a $2$-cell of $X$ labeled by $e\in E$ we can work with the fibration $B\mu_d \to BA(e)\to B\bar A(e)$ instead. Given the standard cell structure on $BA(e)$ it turns out that $f_T^*(\gamma)=-\tau$.  To see this observe that the description of the transgression map implies that the value of $f_T^*(\gamma)$ on the $2$-cell $\Phi^2_e$ can be computed by lifting the boundary of the composite $D^2\xrightarrow{\Phi_e^2} X\xrightarrow{f_T} \bar B(\ZZ/d,G)$, which factors through the subspace $ B\bar A(e)$, to $BA(e)$ and then identifying the homotopy class of the loop in $\pi_1(B\mu_d)\cong\ZZ/d$. The resulting element in $\mu_d$ is precisely $-\tau(e)$ by definition of $f_T$ and the relation $\omega^{-\tau(e)} \prod_{v\in e} T(v)^{\epsilon_e(v)} = I_m$.
}
%
A special case of part (3) is proved in \cite{okay2020homotopical} applicable to hypergraphs with $\epsilon_e(v)=\pm 1$ which has a simply connected topological realization. We sketch an alternative approach for the general case: the class $\gamma_G$ comes from a class in \revtwo{$H^2(B\bar G,\ZZ/d)$, where} $\bar G = G/\mu_d$, still denoted by the same symbol. Let $H\subset G$ denote the discrete subgroup generated by $\set{\rev{T(v)}|\;v\in V}$ together with $\mu_d$. Let $\bar H$ denote the quotient $H/\mu_d$. Since $f_T$ induces the trivial map on $\pi_1$ we can reduce to the case where $\pi_1 (X)=1$ by collapsing the non-contractible loops in $X$.
The composite 
$$X \xrightarrow{f_T} \bar B(\ZZ/d,G) \subset B\bar G$$ factors through a map $X\to B\bar H$.
Since   $\pi_1(X)=1$ and the  homotopy groups of $B\bar H$ vanish above dimension $1$  this map is null homotopic. Therefore using part (2) we have $f_T^*(\gamma_G)=[\tau]=0$.   
\end{proof}

\Ex{\label{ex:Mer}{\rm
Mermin square \cite{Mermin} is the   prominent example of a contextual linear constraint system, i.e., it admits an operator solution but not a scalar solution. Let $P_n$ denote the subgroup in $U(2^n)$ consisting of matrices of the form  $i^{a}A_1\otimes A_2\otimes \cdots\otimes A_n$ where $a\in \ZZ/4$ and each $A_i$ is one  of the Pauli matrices
$$
\rev{I}=\twobytwo{1}{0}{0}{1}\;\;\;\; X=\twobytwo{0}{1}{1}{0}\;\;\;\; Y = \twobytwo{0}{-i}{i}{0} \;\;\;\; Z=\twobytwo{1}{0}{0}{-1}.
$$ 
The linear constraint system $(\hH_{\sq},\tau_{\sq})$ and an operator solution $T_{\sq}:V\to P_2$  is depicted in Figure \ref{fig:Msq} (left figure). 
As depicted in the right figure $\hH_\sq$ has a topological realization given by a torus.
The class $[\tau_{\sq}]$ is non-zero since the cocycle evaluates to $1$ on the torus. Therefore the linear constraint system does not admit a scalar solution \cite{Coho}.

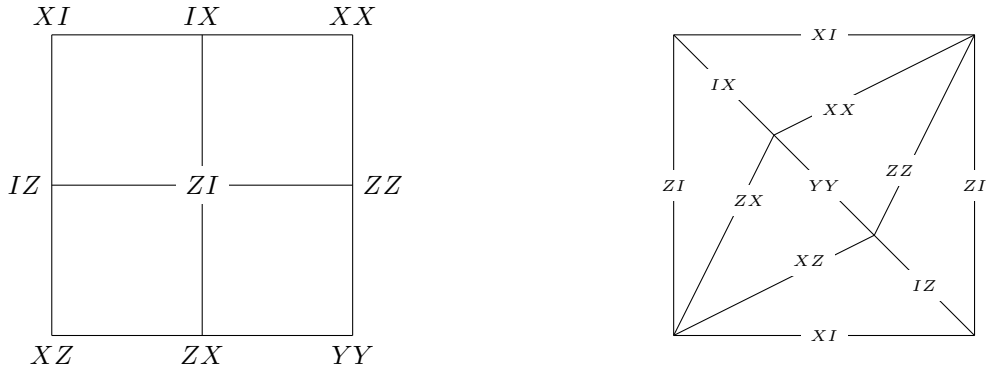
\begin{figure}[h!]
\begin{minipage}{.5\textwidth}
\centering
\begin{tikzpicture}
\draw (0,4) -- (2,4) -- (4,4);
\draw (0,2) -- (2,2) -- (4,2);
\draw (0,0) -- (2,0) -- (4,0);
\draw (0,0) -- (0,2) -- (0,4);
\draw (2,0) -- (2,2) -- (2,4);
\draw (4,0) -- (4,2) -- (4,4);

\node [above] at (0,4) {\footnotesize $XI$};
\node [above] at (2,4) {\footnotesize  $IX$};
\node [above] at (4,4) {\footnotesize  $XX$};

\node [left] at (0,2) {\footnotesize  $IZ$};
\node [fill=white] at (2,2) {\footnotesize  $ZI$};
\node [right] at (4,2) {\footnotesize  $ZZ$};

\node [below] at (0,0) {\footnotesize  $XZ$};
\node [below] at (2,0) {\footnotesize  $ZX$};
\node [below] at (4,0) {\footnotesize  $YY$};

\end{tikzpicture}
\end{minipage}%
\begin{minipage}{.5\textwidth}
\centering
\begin{tikzpicture}
\draw (0,0) -- (4,0);
\draw (4,0) -- (4,4);
\draw (4,4) -- (0,4);
\draw (0,4) -- (0,0);
\draw (0,4) -- (4,0);
\draw (0,0) -- (4/3,8/3); 
\draw (0,0) -- (8/3,4/3);
\draw (4,4) -- (4/3,8/3); 
\draw (4,4) -- (8/3,4/3);

\node [fill=white] at (2,0) {\tiny $XI$};
\node [fill=white] at (2,4) {\tiny $XI$};
\node [fill=white] at (0,2) {\tiny $ZI$};
\node [fill=white] at (4,2) {\tiny $ZI$};
\node [fill=white] at (4/6,5*4/6) {\tiny $IX$};
\node [fill=white] at (5*4/6,4/6) {\tiny $IZ$};
\node [fill=white] at (3*4/6,3*4/6) {\tiny $YY$};
\node [fill=white] at (3*4/12,3*4/12+0.8) {\tiny $ZX$};
\node [fill=white] at (3*4/12+0.8,3*4/12) {\tiny $XZ$};
\node [fill=white] at (3*4/6*3/2-0.8,3*4/6*3/2) {\tiny $XX$};
\node [fill=white] at (3*4/6*3/2,3*4/6*3/2-0.8) {\tiny $ZZ$};

\end{tikzpicture}

\end{minipage}
\caption{(Left figure) $\hH_{\sq}$ consists of $9$ vertices and $6$ edges each consisting of $3$ vertices in each row and  column. \rev{All the incidence weights are equal to $1$.} The operator solution is given by tensor product of two Pauli matrices, where the notation is simplified by omitting $\otimes$. The function $\tau_\sq$ takes the value $0$ for each hyperedge except the right-most column. (Right figure) A topological realization given by a torus together with a cell structure consisting of triangles. The operators are placed on the edges and each triangle corresponds to an hyperedge. The cocycle $\tau_{\sq}$ assigns $0$ to each triangle except $\set{XX,YY,ZZ}$, which is assigned $1$.
} \label{fig:Msq}
\end{figure}

Another linear constraint system constructed in \cite{Mermin} is the Mermin star linear constraint system, which we denote by $(\hH_\st,\tau_\st)$. An operator solution $T_\st:V\to P_3$ is displayed in Figure \ref{fig:Mst} (left figure). The corresponding topological realization is again a torus, but with a different cell structure (right figure); see \cite{Coho}.

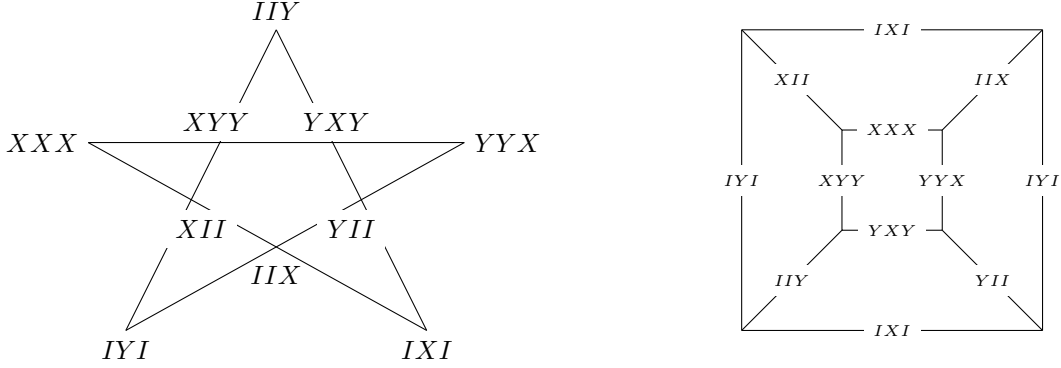
\begin{figure}[h!]
\begin{minipage}{.5\textwidth}
\centering
\begin{tikzpicture}
\draw (0,0) -- (2,4);
\draw (4,0) -- (2,4);
\draw (0,0) -- (4.5,2.5);
\draw (4,0) -- (-0.5,2.5);
\draw (-0.5,2.5) -- (4.5,2.5);
 
\node [below] at (0,0) {\footnotesize  $IYI$};
\node [below] at (4,0) {\footnotesize  $IXI$};

\node [left] at (-0.5,2.5) {\footnotesize  $XXX$};
\node [right] at (4.5,2.5) {\footnotesize  $YYX$};
\node [above,fill=white] at (1.2,2.52) {\footnotesize  $XYY$};
\node [above,fill=white] at (2.8,2.52) {\footnotesize  $YXY$};

\node [below,fill=white ] at (1,1.6) {\footnotesize  $XII$};
\node [below,fill=white ] at (3,1.6) {\footnotesize  $YII$};

\node [above] at (2,4) {\footnotesize  $IIY$};
\node [below] at (2,1) {\footnotesize  $IIX$};

\end{tikzpicture}

\end{minipage}%
\begin{minipage}{.5\textwidth}
\centering
 
\begin{tikzpicture}
\draw (0,0) -- (4,0);
\draw (4,0) -- (4,4);
\draw (4,4) -- (0,4);
\draw (0,4) -- (0,0);
\draw (0,4) -- (4/3,8/3);
\draw (4,0) -- (8/3,4/3);
\draw (4,4) -- (8/3,-4/3+4);
\draw  (0,0) -- (4/3,4/3);
\draw (4/3,4/3) -- (8/3,4/3)-- (8/3,-4/3+4) -- (4/3,8/3) --(4/3,4/3);

\node [fill=white] at (2,0) {\tiny $IXI$};
\node [fill=white] at (2,4) {\tiny $IXI$};
\node [fill=white] at (0,2) {\tiny $IYI$};
\node [fill=white] at (4,2) {\tiny $IYI$};
\node [fill=white] at (4/6,4/6) {\tiny $IIY$};
\node [fill=white] at (4-4/6,4-4/6) {\tiny $IIX$};
\node [fill=white] at (4/6,4-4/6) {\tiny $XII$};
\node [fill=white] at (4-4/6,4/6) {\tiny $YII$};
\node [fill=white] at (4/3,4/2) {\tiny $XYY$};
\node [fill=white] at (4-4/3,4/2) {\tiny $YYX$};
\node [fill=white] at (4/2,4/3) {\tiny $YXY$};
\node [fill=white] at (4-4/2,4-4/3) {\tiny $XXX$};
\end{tikzpicture}
\end{minipage}
\caption{(Left figure) $\hH_{\st}$ consists of $10$ vertices and $5$ edges each consisting of $4$ vertices in each line\rev{, and all the incidence weights are equal to $1$.}  The function $\tau_\st$ takes the value $0$ for each hyperedge except the horizontal line. (Right figure) On the torus $\tau_\st$ specifies a $2$-cocycle that assigns  $0$ to each cell except $\set{XXX,YYX,YXY,XYY}$ is assigned  $1$.  
} \label{fig:Mst}
\end{figure}
}}

\subsection{Computing the homotopy classes}

\Def{Let $X$ be a pointed connected $2$-dimensional CW complex. Consider the collection of triples $(\hH,\tau,T)$ consisting of  a linear constraint system $(\hH,\tau)$ over $\ZZ/d$ where $\hH$ admits a topological realization homotopy equivalent to $X$ and an operator solution $T$ over $G$. Two such triples $(\hH_0,\tau_0,T_0)$ and $(\hH_1,\tau_1,T_1)$ are said to be equivalent if $f_{T_0}$ and $f_{T_1}$ are homotopic as pointed maps. We write $\Sol(X;d,G)$ for the set of equivalence classes and refer to this set as the set of equivalence classes of operator solutions for $(X,d)$ over $G$.
}

The equivalence classes of operator solutions map to the (pointed) homotopy classes of maps
$$
\theta:\Sol(X;d,G)  \hookrightarrow  [X,\bar B(\ZZ/d,G)].
$$
The target can be computed using an algebraic category (the category of {\it crossed modules} \cite{whitehead1949combinatorial}) which captures the behavior of the homotopy category of $2$-dimensional CW complexes.

Let $\bar\pi_i$   denote the $i$-th homotopy group of $\bar B(\ZZ/d,G)$. 

\Pro{\label{pro:class-G}
Let $X$ be a connected $2$-dimensional CW complex. 
Sending a map to the homomorphism induced on $\pi_1$ gives a surjective map
$$
\pi:[X,\bar B(\ZZ/d,G)] \to \Hom(\pi_1X,\bar\pi_1)
$$
such that  for a fixed homomorphism $\alpha$ the preimage is given by
$$
\pi^{-1}(\alpha) \cong H^2(\tilde X,(\bar\pi_2)_\alpha)
$$
where $(\bar\pi_2)_\alpha$ is the $\pi_1(X)$-module determined by the homomorphism $\alpha$. 
}  
\Proof{
The statement holds for $[X,Y]$ where $Y$ is an arbitrary CW complex. We will construct maps
$$
Y \xrightarrow{r} \bar Y \xleftarrow{s} Y_{(2)}
$$
where $Y_{(2)}$ is a $2$-dimensional CW complex, and the maps $r$ and $s$ are  $3$-equivalences, i.e., each map  induces an isomorphism on $\pi_i$ for $0\leq i<3$ and a surjection for $i=3$. In this case $r_*:[X,Y] \to [X,\bar Y]$, and similarly $s_*$,  are bijections \rev{\cite[Cor. 23, p. 405]{spanier1989algebraic}}. Before the construction we first show how to finish the proof of the statement.

The set $[X,Y_{(2)}]$   can be computed algebraically; for details we refer to \rev{\cite[Ch. II, \S 4.2]{TwoDim}}. Let us write $[X,Y_{(2)}]_\alpha$ for the set of  homotopy classes of maps that induce the homomorphism $\alpha$ between the fundamental groups.  The (cellular) chain complex for the universal cover $\tilde X$ consists of $\pi_1(X)$-modules and we can talk about the cohomology groups $H^n(\tilde X, (\pi_2 Y_{(2)})_\alpha)$ where $\pi_2 Y_{(2)}$ is regarded as a $\pi_1(X)$-module via the homomorphism $\alpha$.  The cohomology group $H^2(\tilde X, (\pi_2   Y_{(2)})_\alpha)$ acts on $[X,  Y_{(2)}]_\alpha$ in a transitive way \rev{\cite[Ch. II, Thm. 4.11]{TwoDim}) and} this action determines a bijection  
$$
[X,  Y_{(2)}]_\alpha \cong H^2(\tilde X, (\pi_2 Y_{(2)})_\alpha).
$$

We turn to the construction of $r$ and $s$. The first map is obtained by killing homotopy groups of $Y$ above dimension $2$. Construction of the second map uses the theory of crossed modules.  The fundamental property   we will use  is that any free crossed module over a free base group is realizable by a $2$-dimensional CW complex and maps between such crossed modules come from maps between the CW complexes that realize them \cite[Ch. II]{TwoDim}. Let us apply this to the crossed module given by the connecting homomorphism
\begin{equation}\label{eq:crossed}
\partial:\pi_2(\bar Y,\bar Y^1) \to \pi_1(\bar Y^1)
\end{equation}
By the realization result there is a $2$-dimensional CW complex $Y_{(2)}$ such that the crossed module  $\partial:\pi_2(Y_{(2)},Y^1_{(2)}) \to \pi_1(Y^1_{(2)})$ is isomorphic to the one given in \ref{eq:crossed}.  
We will show that this isomorphism is realized by a map 
$
s: Y_{(2)} \to \bar Y.
$ 
We start the construction of $s$ from the $1$-st skeleton. We can find a map $Y_{(2)}^1 \to \bar Y^1$ that induces the desired isomorphism on $\pi_1$. Composing this map with the inclusion $\bar Y^1\subset \bar Y$ we obtain $Y_{(2)}^1\to \bar Y$.  This map lifts to a map $Y_{(2)}\to \bar Y$ since the set of $2$-cells is a basis for the free group $\pi_2(Y_{(2)},Y_{(2)}^1)$ and the isomorphism between the crossed modules implies precisely the lifting condition in the algebraic language. 
}

\begin{ex}{\rm  \label{ex:pauli}
\rev{We will discuss an interesting example related to the Pauli group. For properties of this group we refer to \cite{okay2019classifying}.}
The Pauli group $P_n$ defined in Example \ref{ex:Mer} has a generalization for all primes $p$ which has a similar description as tensor products of $p\times p$ unitary matrices. As an abstract group $P_n$ is \rev{an almost extraspecial $2$-group for $p=2$} 
and an extraspecial $p$-group of exponent $p$ for odd primes. 
There is an irreducible complex representation \rev{of $P_n$} which allows us to regard it as a subgroup in $U(p^n)$.

\rev{The central quotient group of $P_n$, which will be denoted by $E_n$, is an elementary abelian $p$-group of rank $2n$. There is a symplectic bilinear form $\bi$ on $E_n$ induced by the commutator of $P_n$. We can choose a symplectic basis $\set{x_1,x_2,\cdots,x_n,z_1,z_2,\cdots,z_n}$ for $E_n$. Let $B(\bi,E_n)$ denote the geometric realization of the simplicial set
$$
[k]\mapsto \set{(a_1,a_2,\cdots,a_k)\in (E_n)^k |\;\bi(a_i,a_j)=0\;\forall i,j }.
$$
The simplicial structure is induced from the simplicial set $E_n^\dt$ whose geometric realization is the classifying space $BE_n$. The space $\bar B(\ZZ/p,P_n)$ can be identified with $B(\bi,E_n)$.
}
\rev{For $n\geq 2$ i}t is known that 
$$
\pi_1 \bar B(\ZZ/p,P_n) \rev{= \pi_1 B(\bi,E_n) } \cong \left\lbrace \begin{array}{ll}
\ZZ/2 \times \rev{E_n} & p=2 \\
P_n  & p>2,
\end{array} \right.
$$
and the higher homotopy groups are given by 
$$
\pi_i \bar B(\ZZ/p,P_n)\rev{= \pi_i B(\bi,E_n) } \cong \pi_i(\bigvee^{N_{p,n}}S^n),\;\;\; i\geq 2,
$$ 
where $N_{p,n}$ has an explicit  formula \rev{\cite[\S 6]{O16}}.  Therefore according to Proposition \ref{pro:class-G} the map
$$
[X,\bar B(\ZZ/p,P_n)] \to \Hom(\pi_1 X , \bar\pi_1)
$$
is an isomorphism when $n\geq 3$. However, for $n=2$ 
it is only surjective and the kernel depends on the $\bar\pi_1$-module structure of $\bar\pi_2$, which is currently unknown.

The canonical class can be described as
\begin{equation}\label{eq:cup}
\gamma_{P_n} = \left\lbrace
\begin{array}{cc}
 x_0^2 + \sum_{i=1}^n x_i\cup z_i & p=2\\
0 & p>2,
\end{array}
\right.
\end{equation}
where $\set{x_0,x_1,\cdots,x_n,\rev{z_1},\cdots,z_n}$ is a basis for $\ZZ/2 \times \rev{E_n}$; see \cite{okay2019classifying} for details. Therefore for odd $p$ every linear constraint system has a scalar solution if it has an operator solution over $P_n$. Whereas  for $p=2$ this depends on the map induced on $\pi_1$, as a result of the cup product decomposition in \ref{eq:cup} .

The operator solution $T_{\sq}$ of the Mermin square linear system $(\hH_{\sq},\tau_{\sq})$ introduced in Example \ref{ex:Mer} gives a non-trivial class $[f_{T_{\sq}}]$ in $[S^1\times S^1, \bar B(\ZZ/2,P_2)]$. For $n\geq 2$ let us write 
\begin{equation}\label{eq:Tn}
T_{n}=T_{\sq}\otimes I_{2^{n-2}},
\end{equation}
 for the operator solution obtained by tensoring with the identity matrix: $A\mapsto A\otimes \rev{I_{2^{n-2}}}$. Then $[f_{T_n}]$ gives a non-trivial class in $[S^1\times S^1, \bar B(\ZZ/2,P_n)]$ for all $n\geq 2$.
Similarly the Mermin star example $(\hH_\st,\tau_\st)$ specifies a class in $[S^1\times S^1,\bar B(\ZZ/2,P_3)]$. It turns out that this class coincides with $[f_{T_3}]$ since there is a refined cell structure (\cite{Coho}) on the torus as depicted in Figure \ref{fig:refined}.
More precisely, there is a commutative diagram
$$
\begin{tikzcd}
X_{sq} \arrow{d} \arrow[r,hook] & X_{ref}\arrow{d}  &\arrow[l,hook']  X_{st} \arrow{d} \\  
 \bar B(\ZZ/2,P_2) \arrow{r}{\otimes I_2}& \bar B(\ZZ/2,P_3)  &\arrow[l,equal] \bar B(\ZZ/2,P_3)
\end{tikzcd}
$$
relating the topological realizations $X=S^1\times S^1$ with different cell structures as indicated by the subscripts.

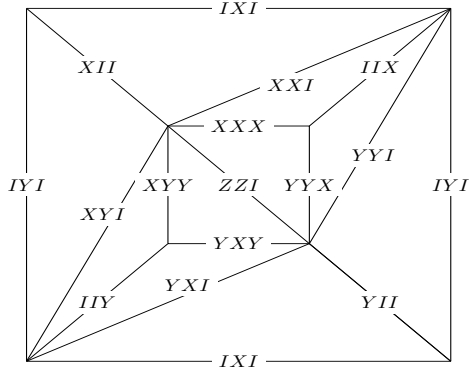
\begin{figure}[h!]

\begin{center}
\resizebox{6.5cm}{5cm}{
\begin{tikzpicture} 
\draw (2*0,2*0) -- (2*4,2*0);
\draw (2*4,2*0) -- (2*4,2*4);
\draw (2*4,2*4) -- (2*0,2*4);
\draw (2*0,2*4) -- (2*0,2*0);
\draw (2*0,2*4) -- (2*4/3,2*8/3);
\draw (2*4,0) -- (2*8/3,2*4/3);
\draw (2*4/3,2*8/3) -- (2*4,0);
\draw (2*4,2*4) -- (2*8/3,-2*4/3+2*4);
\draw  (0,0) -- (2*4/3,2*4/3);
\draw (2*4/3,2*4/3) -- (2*8/3,2*4/3)-- (2*8/3,-2*4/3+2*4) -- (2*4/3,2*8/3) --(2*4/3,2*4/3);
\draw (0,0) -- (2*4/3,2*8/3); 
\draw (0,0) -- (2*8/3,2*4/3);
\draw (2*4,2*4) -- (2*4/3,2*8/3); 
\draw (2*4,2*4) -- (2*8/3,2*4/3);

\node [fill=white] at (2*2,0) {$IXI$};
\node [fill=white] at (2*2,2*4) {$IXI$};
\node [fill=white] at (0,2*2) {$IYI$};
\node [fill=white] at (2*4,2*2) {$IYI$};
\node [fill=white] at (2*4/6,2*4/6) {$IIY$};
\node [fill=white] at (2*4-2*4/6,2*4-2*4/6) {$IIX$};
\node [fill=white] at (2*4/6,2*4-2*4/6) {$XII$};
\node [fill=white] at (1.5*16/6,1.5*16/6) {$ZZI$};
\node [fill=white] at (2*4-2*4/6,2*4/6) {$YII$};
\node [fill=white] at (2*4/3,2*4/2) {$XYY$};
\node [fill=white] at (2*4-2*4/3,2*4/2) {$YYX$};
\node [fill=white] at (2*4/2,2*4/3) {$YXY$};
\node [fill=white] at (2*4-2*4/2,2*4-2*4/3) {$XXX$};
\node [fill=white, below left] at (2*3*4/12,2*3*4/12+2*0.8) {$XYI$};
\node [fill=white, below left] at (2*3*4/12+2*0.8,2*3*4/12) {$YXI$};
\node [fill=white,above right] at (2*3*4/6*3/2-2*0.8,2*3*4/6*3/2) {$XXI$};
\node [fill=white,above right] at (2*3*4/6*3/2,2*3*4/6*3/2-2*0.8) {$YYI$};
\end{tikzpicture}
}
\end{center}

\caption{Refined topological realization} \label{fig:refined}
\end{figure}  

}
\end{ex}

\subsection{Application of $C(d,m)$-cohomology} Now we focus on operator solutions in unitary groups. For notational simplicity let us write $\Sol(X;d,m)$ for \rev{the set of equivalence classes of operator solutions} over $U(m)$. Recall the map 
$$\bar\iota_m:\bar B(\ZZ/d,U(m)) \to \bar B(d,m)$$
 introduced in \ref{eq:diag-Q-fib}. Composing with $\bar\iota_m$ gives a map
$$
\hat\theta:\Sol(X;d,m) \hookrightarrow [X,\bar B(\ZZ/d,U(m))] \xrightarrow{(\bar\iota_m)_*} C(d,m)(X)
$$
where we have identified $[X,\bar B(d,m)]$ with the $0$-th $C(d,m)$-cohomology of $X$ since the target space is the infinite loop space associated to the spectrum representing the cohomology theory. Given an operator solution $T$ the image of $f_T$ under $(\bar\iota_m)_*$ will be denoted by $\hat f_T$. 

\Rem{{\rm By Lemma \ref{lem:stablecan} the pull-back $\hat f_T^*(\gamma_m^\sS)$ coincides with $f_T^*(\gamma_m)$. Therefore for a linear constraint system existence of a scalar solution is determined in a stable manner, i.e., $\hat f_T^*(\gamma_m^\sS)=0$ if and only if a scalar solution exists.
}}


\Cor{\label{cor:stable-LCS}
Let $(\hH,\tau)$ be a linear constraint system  over $\ZZ/d$ and $X$ be a topological realization for $\hH$.
\begin{enumerate}
\item If \rev{$(\hH,\tau)$ has an operator solution and} $H^2(X,(\ZZ/d)_m)=0$ then  $(\hH,\tau)$   has a scalar solution.
\item If $d$ and $m$ are coprime then $C(d,m)(X)=0$. In particular,  $(\hH,\tau)$    has a scalar solution if it has an operator solution over $U(m)$.
\item If $\pi_1(X)$ is trivial and $[\tau]\not= 0$ then $(\hH,\tau)$ does not have an operator solution over $U(m)$ for any \rev{$m\geq 2$}. 

\end{enumerate}
}
\begin{proof}\rev{Let $T$ be an operator solution for $(\hH,\tau)$ and $f:\sS\wedge X\to C(d,m)$ represent the class $\hat\theta(T)$.} If $H^2(X, (\ZZ/d)_m)=0$ then by \revtwo{Lemma \ref{lem:stablecan} part (2) $f_T^*(\gamma_m)=0$ and the fibration $B\mu_d \to \tilde X\to X$ classified by $f_T^*(\gamma_m)$ splits, i.e $\tilde X\simeq X \times B\mu_d$. Choosing a splitting $X\to \tilde X$ and composing with the map $\tilde X\to B(\ZZ/d,U) \xrightarrow{\det} B\mu_d$ shows that there is a class $\alpha$ in $k\mu_d(X)$ such that $[f]=(\zeta(\alpha),\cl(f))$, where $\cl(f)$ is regarded as an element of $H^2(X,(\ZZ/d)_m)$, under the decomposition given in Theorem \ref{thm:Thm-stable}. Since $H^2(X,(\ZZ/d)_m)=0$ we have $\cl(f)=0$ and}
(1) follows from 
Proposition \ref{pro:LCS-basics} part (1) and (2).  
Part (2) follows from part (1) since if $(d,m)=1$ then $H^2(X,(\ZZ/d)_m)=0$. 
Part (3) follows from Proposition \ref{pro:LCS-basics} part (3).
  Existence of an operator solution implies that $[\tau]=0$ since $X$ is simply connected.
\end{proof}

\subsection{The Mermin class}\label{sec:Mermin-class}
 In the physics literature a quantum system with Hilbert space  $ (\CC^2)^{\otimes n}$ is called an {\it $n$-qubit}. Such systems play a significant role in quantum information theory. Operator solutions in $U(2^n)$ of linear constraint systems over $\ZZ/2$ produce classes in $C(d,m)$-cohomology, where $d=2$ and $m=2^n$. Theorem \ref{thm:Thm-stable} gives an isomorphism
\begin{equation}\label{eq:C22n-iso}
C(2,2^n)(X) \cong H^1(X,\ZZ/2)\oplus H^2(X,\ZZ/2).
\end{equation}
We will construct non-trivial classes that come from operator solutions of the Mermin square linear constraint system described in Example \ref{ex:Mer}. Our topological realization is a torus $X=S^1\times S^1$. An operator solution for $n=2$ is given in Figure \ref{fig:Msq}. Let $T_1$ denote this solution. We define an operator solution in $U(2^n)$ by tensoring with the identity as in \ref{eq:Tn}, i.e., by constructing an operator solution $T_n$ defined by $T_n(v)=T_1(v)\otimes I_{2^{n-1}}$ for $v\in V$. Let $[T_n]$ denote the class of this solution in $\Sol(S^1\times S^1;2,2^n)$. Let $M_n$ denote the class $\hat{\theta}(T_n)$ in $C(2,2^n)(S^1\times S^1)$. This class will be called the {\it Mermin class}. \rev{As a consequence of the isomorphism \ref{eq:C22n-iso} $M_n$ is represented by a pair of cohomology classes} 
$$\rev{(\varphi_1;\varphi_2)\in  H^1(X,\ZZ/2)\oplus H^2(X,\ZZ/2).}$$ 
For each $n\geq 2$ the cohomology class $[\tau]\neq 0$ since, as we have seen in Example \ref{ex:Mer}, the Mermin square linear constraint system does not admit a scalar solution. Therefore $\varphi_2$ is the non-trivial class in $H^2(S^1\times S^1,\ZZ/2)=\ZZ/2$. 
To determine $\varphi_1$  
\rev{we will construct a homotopy commutative}  diagram  
\begin{equation}\label{eq:Xphi2}
\begin{tikzcd}
X_{\varphi_2} \arrow{r} \arrow{d} & B(\ZZ/2,U) \arrow{d} \\
X\arrow{r} & \bar B(2,2^n)
\end{tikzcd}
\end{equation}
where  $X_{\varphi_2}$ is the $B\mu_2$-bundle determined by the non-trivial class $\varphi_2$. \rev{The bottom map in the diagram is given by the composite
\begin{equation}\label{eq:long-composite}
X\xrightarrow{f_{T_n}} \bar B(\ZZ/2,P_n) \to \bar B(\ZZ/2, U(2^n)) \xrightarrow{\bar \iota_{2^n}} \bar B(2,2^n) 
\end{equation}
where the middle map is induced by the inclusion $P_n\subset U(2^n)$. Note that the map obtained using Construction \ref{con:fT} factors through $\bar B(\ZZ/2,P_n)$ since the operator solution $T_n$ is over the Pauli group $P_n$. Each map in \ref{eq:long-composite} can be extended to a map between homotopy fibrations over $B^2\mu_2$ corresponding to the cohomology classes $\varphi_2$, $\gamma_{P_n}$ (see \ref{eq:cup}), $\gamma_{2^n}$ and $\gamma^\sS_{2^n}$; respectively. The top map in \ref{eq:Xphi2} is given by the composite map induced between the homotopy fibers. Now we observe that t}he class $\varphi_1$ is determined by the map induced on $\pi_1X\to \pi_1 \bar B(2,2^n)$. Applying $\pi_1$ to the diagram in \rev{\ref{eq:Xphi2}} we obtain
$$
\begin{tikzcd}
\pi_1 X_{\varphi_2} \arrow{r} \arrow[two heads,d] & \pi_1 B(\ZZ/2,U) \arrow[d,"\cong"] \\
\pi_1X\arrow{r} & \pi_1\bar B(2,2^n)
\end{tikzcd}
$$   
Let $\tilde x$ and $\tilde z$ denote the elements lifting the generators $x=(1,0)$ and $z=(0,1)$ of the quotient group $\pi_1 X=\ZZ^2$. It suffices to determine the images of $\tilde x$ and $\tilde z$ under the top horizontal map. Figure \ref{fig:Msq} tells us that $\tilde x$ maps to the loop determined by $X\otimes I_{2^{n-1}}$ and $\tilde z$ maps to $Z\otimes I_{2^{n-1}}$. We can understand the induced map on $\pi_1$ by composing with the determinant map $\det:B(\ZZ/2,U)\to B\mu_2$. This amounts to taking the determinant of the matrices representing the loops, which gives $1$ in both cases. Thus both of the loops map to the trivial loop in $B\mu_2$. 
Therefore $\varphi_1=(0,0)\in   (\ZZ/2)^2$ \revtwo{(also follows from Corollary \ref{cor:iota-bar-alternative})}. In summary, the Mermin class $M_n$ is represented by $(0,0;1)$.    
Since $f_{T_n}$ induces the trivial map on $\pi_1$ it factors as
\begin{equation}\label{eq:S2fac}
\begin{tikzcd}
S^1\times S^1 \arrow{d} \arrow{r}{f_{T_n}} & \bar B(2,2^n)\\
S^2 \arrow[dashed,ru,"\bar f"']
\end{tikzcd}
\end{equation}
where the vertical map collapses the non-contractible loops corresponding to $x$ and $z$. The homotopy class of $\bar f$ is the generator of $\pi_2 C(2,2^n)=\ZZ/2$.
By slight abuse of notation we will also write $M_n$ for this class and   refer to it as the Mermin class as well.


Let us compare to  the unstable situation. 
\rev{By looking at the homotopy fibers of the maps in \ref{eq:long-composite} regarded as homotopy fibrations over $B\mu_2$ as before we see that t}he diagram \ref{eq:Xphi2} factors as
$$
\begin{tikzcd}
X_{\varphi_2} \arrow{r} \arrow{d} & B(\ZZ/2,P_n) \arrow{r}{\tilde g} \arrow{d}  & B(\ZZ/2,U) \arrow{d}   \\
X\arrow{r}{f} & \bar B(\ZZ/2,P_n) \arrow{r}{g} & \bar B(2,2^n)  
\end{tikzcd}
$$
where $g^*(\gamma_{2^n}^\sS) = \gamma_{P_n}$. The homotopy class $[f]$ is non-trivial in $[X,\bar B(\ZZ/2,P_n)]$, which surjects onto $\Hom(\pi_1 X, \bar \pi_1)$ as we have seen in Example \ref{ex:pauli}. However, the composite $gf$ induces the trivial map on $\pi_1$. This is not in conflict with Proposition \ref{pro:LCS-basics} part (3) if we take $G=U$. This is because the subgroup $\mu_2\hookrightarrow U$ is not a central, or even not a  normal, subgroup. Proposition \ref{pro:LCS-basics} part (3) also implies that the diagonal map $\bar f$ in \ref{eq:S2fac} does not factor through  $\bar B(\ZZ/2, U(2^n))$.

\subsection{Relation to \rev{SPT} phases}\label{sec:ex-real} 
\rev{Replacing $U(m)$ with the orthogonal group $O(m)$
gives the real versions of the spectra considered in this paper. As explained in Appendix \ref{app:real-case}   we can define the real symmetric $K$-theory spectrum $ko_\sym$ and the spectrum $C_\RR(m)$, the real version of the $C(d,m)$ spectrum.}
Mermin square construction and its $n$-qubit version $T_n$ \rev{introduced in \S\ref{sec:Mermin-class}} can be regarded as an operator solution over $O(2^n)$ since the matrices involved have real entries.  Let $M_n^\RR\in C_\RR(2^n)(S^1\times S^1)$ denote the corresponding class. 
As in the complex case we find that $M_n^\RR$ can be identified with the generator of the quotient in the exact sequence \rev{(see \ref{eq:exact-homotopy-Cdm-real})}
$$
0\to \pi_2 ko_\sym  \to  \pi_2 C_\RR(2^n) \to H^2(S^2,\ZZ/2) \to  0.
$$
\rev{T}he generator of $\pi_2 ko_\sym=\ZZ/2$ has also a physical interpretation. 
It can be realized as a non-trivial SPT phase: The $ko$-orientation $M\Spin\to ko$ of the  spin cobordism spectrum $M\Spin$ is highly connected \rev{(\revtwo{this} follows from \cite[Thm. 2.2]{anderson1967structure}; see also \cite[Rem. 6.1]{campbell2017homotopy})}. In particular, it induces an isomorphism on $\pi_2$.
Therefore smashing this map with $B\mu_2$ induces an isomorphism $\pi_2(M\Spin\wedge B\mu_2) \to \pi_2(ko\wedge B\mu_2)$. The generator of $\pi_2(M\Spin\wedge B\mu_2)$ is identified as the Gu--Wen phase, a fermionic SPT phase constructed in \cite{gu2014symmetry}; see also \cite[\S 5]{kapustin2015fermionic}.   This class hits the generator of $\pi_2 ko_{\sym}$ under the identification $\pi_2(ko\wedge B\mu_2) \cong \pi_2(ko_\sym)$.

\appendix

\renewcommand{\theequation}{\thesection.\arabic{equation}}

\section{Real version}\label{app:real-case}

\rev{

In this section we describe the real versions of the constructions introduced in the previous sections. The main idea is to replace $U(m)$ with the orthogonal group $O(m)$. As in the complex case we obtain a commutative version of the topological real $K$-theory and a cohomology theory, denoted by $C_\RR(m)$, that can be used to study operator solutions of LCSs over $O(m)$. 
}

Every abelian subgroup of $O(m)$ can be conjugated into $\SO(2)^j \times O(1)^{m-2j}$ \rev{for some $j\leq  \lfloor  m/2 \rfloor$} \cite[Appendix A]{higuera2014spaces}. Thus a homomorphism $f:\ZZ^m\to O(m)$, when regarded as a representation, is isomorphic to a direct sum
$$
f\cong \eta_1\oplus \eta_2\oplus \cdots \oplus \eta_j \oplus \ell_1 \oplus \ell_2\oplus \cdots \oplus \ell_{2m-j}
$$
where $\eta_i:\ZZ^m\to SO(2)$ and $\ell_i:\ZZ^m\to O(1)$. In particular, a matrix is diagonalizable in $O(m)$ if and only if it is {\it symmetric}, i.e. $A^T=A$. Thus in the real case we will consider $2$-torsion orthogonal matrices. The resulting space $B(\ZZ/2,O(m))$ is constructed from pairwise commuting symmetric orthogonal matrices. We can stabilize over $m$, similar to the complex case, to obtain $B(\ZZ/2,O)$. \rev{\revtwo{The} real version of  Proposition \ref{pro:B-ku} gives a homeomorphism
$$
ko((\mu_2)^n) \xrightarrow{\cong} \Hom((\ZZ/2)^n,O)
$$
where $ko$ is the corresponding $\Gamma$-space of the connective real $K$-theory spectrum.
Similar to the complex case this homeomorphism  is compatible with the simplicial structures and induces a homeomorphism
$$
ko(B\mu_2) \xrightarrow{\cong} B(\ZZ/2,O).
$$  
From $ko$ we can construct the $\Gamma$-space $ko_{B\mu_2}$ and consider t}he associated spectrum  
$$
ko_{\sym} = ko_{B\mu_2}(\sS).
$$ 
We will refer to $ko_\sym$ as the {\it real symmetric $K$-theory}. There is a similar stable equivalence $ko_\sym\simeq ko\wedge B\mu_2$ and a weak equivalence $B(\ZZ/2,O)\simeq \Omega^\infty ko_\sym$ by the real versions of Proposition \ref{pro:B-ku} and \ref{pro:InfLoop} \rev{(see also \cite[Rem. 2.9]{gritschacher2019commuting})}. 
\revtwo{The} homotopy groups of $B(\ZZ/2,O)$ are isomorphic  to the \rev{(reduced)} $ko$-homology of $B\mu_2$:
\begin{equation}\label{tab:ko-BZ/2}
\pi_{8k+\epsilon}(ko_\sym) = \left\lbrace 
\begin{array}{cc}
\ZZ/2  & \epsilon = 1,2 \\
\ZZ/2^{4k+3} & \epsilon =3 \\
\ZZ/2^{4k+4} & \epsilon =7 \\
0 & \text{otherwise.}
\end{array}
\right.
\end{equation}
\rev{(taken from the unreduced version in \cite[\S 12.2.D]{bruner2010connective}).}
 
Similar to the complex case $\pi_1(ko_\sym)$ can be understood by considering the composition of $B\mu_2 \subset B(\ZZ/2,O)$ with the determinant map $\Det:B(\ZZ/2,O)\to B\mu_2$. This composition is the identity map and splits off the $\ZZ/2$ in the first homotopy group. Moreover, the unit map $\sS \to ko$ is $3$-connected, i.e. induces an isomorphism on $\pi_i$ for \rev{$0\leq i< 3$} and a surjection on $i=3$\rev{, (mainly because $\pi_1(ko)$ is generated by the image of the Hopf map $\eta\in \pi_1(\sS)$ \cite[Thm. 3.1.26]{ravenel2003complex}).}  From the Atiyah--Hirzebruch spectral sequence we see that $\sS\wedge B\mu_2 \to ko\wedge B\mu_2$ is also $3$-connected. \rev{In addition, since $\pi_3(QB\mu_2)\cong \ZZ/8$ \cite{liulevicius1963theorem}} the map 
$
 Q(B\mu_2) \to B(\ZZ/2,O)
$
extending the inclusion $B\mu_2 \subset B(\ZZ/2,O)$ induces an isomorphism on $\pi_r$ for $0\leq r\leq 3$. \rev{Therefore we have
\begin{equation}\label{eq:QBmu2-homotopy}
\pi_r^s(B\mu_2)= \pi_r(QB\mu_2) \cong \left\lbrace \begin{array}{cc}
0 & r=0 \\
\ZZ/2 & r=1,2 \\
\ZZ/8 & r=3. \\
\end{array} \right.
\end{equation}
}

\rev{ 
Let $\delta_m:\sS\to ko$ denote the morphism of $\Gamma$-spaces corresponding to the real version of \ref{eq:gamma-map}. There is an induced map of spectra $\delta_m\wedge \idy:\sS\wedge B\mu_2 \to ko\wedge B\mu_2$.  
 
\Def{\label{def:real-Cm}{\rm
The spectrum \rev{$C_\RR(m)$} is obtained by killing the homotopy groups above degree $2$ of  the cofiber of $\delta_m\wedge \idy$.
}} 
 
} 

Theorem \ref{thm:Thm-stable} has also a real version  where \rev{$C(d,m)$} is replaced by \rev{$C_\RR(m)$}.

\rev{
\Lem{\label{lem:real-multby-m} Let  $(\delta_m\wedge \idy)_*:\pi_n(\sS\wedge B\mu_2) \to \pi_n(ko\wedge B\mu_2)$ denote the homomorphism between the homotopy groups induced by the spectrum map $\delta_m\wedge \idy:\sS\wedge B\mu_2 \to ko\wedge B\mu_2$. Then
$$
(\delta_m\wedge  \idy)_* =  \underbrace{(\delta_1\wedge \idy)_*+\cdots+(\delta_1\wedge \idy)_*}_m.
$$
} }
\begin{proof}
We \rev{first show} that   $\delta_m:\sS \to \rev{ko}$ is the $m$-fold sum $\delta_1+\cdots+\delta_1$.  By definition \rev{the $\Gamma$-space morphism} $\delta_1$ is completely determined by its value on $\sS(1_+)=1_+$, which sends $1$ to the subspace $\Span{e_1}\rev{\subset\RR^\infty}$. We have an $H$-space structure on $\rev{ko}(1_+)$, which comes from being a special $\Gamma$-space, that is induced by 
\begin{equation}\label{eq:direct-sum}
\rev{ko}(1_+) \times \rev{ko}(1_+) \to \rev{ko}(1_+)
\end{equation}
that sends $(V,W)$ to the direct sum $V\oplus W$. This $H$-space structure is responsible for the abelian group structure on the set of homotopy classes of maps $[\sS,\rev{ko}]$. Thus $\delta_1+\delta_1$ is computed by using \ref{eq:direct-sum}. In effect we obtain a map $\sS(1_+)\to \rev{ko}(1_+)$ that sends $1$ to the direct sum $\Span{e_1}\oplus \Span{e_1}\cong \Span{e_1,e_2}$. This is precisely $\delta_2$. In a similar way we can proceed to show that $\delta_m$ is the $m$-fold sum of $\delta_1$  as claimed. 
\rev{To finish the proof} we rely on the following basic properties of the homotopy category of spectra: 
Let $K,L,M$ be spectra, $X$ be a space, and $f,f':L\to M  $ be  maps of spectra. 
\begin{enumerate}

\item $\wedge\idy: [K,L] \to [K\wedge X, L\wedge X]$, defined by $f\mapsto f\wedge \idy$, is a homomorphism of abelian groups, i.e, $(f+f')\wedge \idy= f\wedge \idy + f'\wedge \idy$.  

\item Consider the induced map $f_*:[K,L]\to [K,M]$, defined by $f_*(g)=f g$, and $f'_*$ similarly defined. Then $(f+f')_* = f_* + f'_*$.

\end{enumerate}
Both of these results follow from the basic properties of addition of spectrum maps.
We apply (1) to $[\sS,\rev{ko}]\to [\sS\wedge B\rev{\mu_2},\rev{ko}\wedge B\rev{\mu_2}]$ and obtain
\begin{equation}\label{eq:dec}
\delta_m\wedge \idy=(\delta_1+\cdots+\delta_1)\wedge \idy = (\delta_1\wedge \idy)+\cdots +(\delta_1\wedge\idy).
\end{equation}
Note that the map induced on $\rev{\pi_n}$ can be thought of as a map
\begin{equation}\label{eq:Sigma1}
(\delta_m\wedge \idy)_*:[\rev{\Sigma^n}\sS,\sS\wedge B\rev{\mu_2}] \to [\rev{\Sigma^n}\sS,\rev{ko}\wedge B\rev{\mu_2}].
\end{equation} 
Now we apply (2) to the decomposition given in \ref{eq:dec} \rev{and} obtain that
$$
(\delta_m\wedge \idy)_*= (\delta_1\wedge \idy)_*+\cdots +(\delta_1\wedge\idy)_*.
$$
\end{proof}  

\rev{A similar result applies to the complex version and can be used to give an alternative proof of Lemma \ref{lem:Cdm-homotopy}. In the real case  $\sS\wedge B\mu_2 \to ko\wedge B\mu_2$ induces an isomorphism on $\pi_r$ for $0\leq r\leq 3$ as observed above. Then using Lemma \ref{lem:real-multby-m} together with} the homotopy groups of $ko_\sym$ given in \ref{tab:ko-BZ/2} and of $Q(B\mu_2)$ \rev{given in \ref{eq:QBmu2-homotopy} we obtain an exact sequence
$$
0\to \ZZ/2 \xrightarrow{\times m} \ZZ/2 \xrightarrow{\alpha} \pi_2 C_\RR(m) \xrightarrow{\beta}  \ZZ/2 \rev{\xrightarrow{\times m}}  \ZZ/2 \to \pi_1 C_\RR(m) \to 0  .
$$ 
W}e see that if $m$ is odd then $\pi_i \rev{C_\RR(m)}=0$ for $i=1,2$. Thus the interesting case is \rev{when $m$ is even.} 
\rev{In this case t}he homotopy groups fit into the exact sequence
\begin{equation}\label{eq:exact-homotopy-Cdm-real}
0\to \ZZ/2 \xrightarrow{\alpha} \pi_2 C_\RR(m) \xrightarrow{\beta}  \ZZ/2 \xrightarrow{0}  \ZZ/2 \rev{\xrightarrow{\cong}} \pi_1 C_\RR(m) \to 0  .
\end{equation}
Let $X$ be a connected $2$-dimensional CW complex. \rev{For $m$ even t}here is a commutative diagram 
$$
\begin{tikzcd}
H^1(X,\ZZ/2) \arrow{r}{\delta} &H^2(X,\ZZ/2) \arrow[d,hook] \arrow[r, "\alpha_*"] & H^2(X,\pi_2 \rev{C_\RR(m)}) \arrow[d,hook] \arrow[r,"\beta_*",two heads] & H^2(X,\ZZ/2) \arrow[d,equal] \\
&ko_\sym(X) \arrow[d,two heads] \arrow[r,"\zeta"] & \rev{C_\RR(m)}(X) \arrow[d,two heads] \arrow[r,"\cl"] & H^2(X,\ZZ/2) \\ 
&H^1(X,\ZZ/2) \arrow[r,"\cong"] & H^1(X,\ZZ/2) &
\end{tikzcd}
$$
where $\delta$ is the connecting homomorphism of the exact sequence associated to $0\to \ZZ/2 \to \pi_2 C_\RR(m)\to \ZZ/2\to 0$. The image of $\zeta$ is contained in the kernel of $\cl$. 

\begin{acknowledgments} 
The author would like to thank
Simon  Gritschacher for his comments on an earlier version of this paper; Daniel Sheinbaum for discussions on symmetry-protected topological phases\revtwo{; and the anonymous referee for pointing out the splitting in Corollary \ref{cor:Cdm-cofiber} and for Corollary \ref{cor:iota-bar-alternative}.}
This work is supported by \rev{the Air Force Office of Scientific Research under award
number FA9550-21-1-0002.}
\end{acknowledgments} 

Data sharing is not applicable to this article as no new data were created or analyzed in this study.


\begin{thebibliography}{10}

\bibitem{AGLT}
A.~Adem, J.~G{\'o}mez, J.~Lind, and U.~Tillmann.
\newblock Infinite loop spaces and nilpotent {$K$}--theory.
\newblock {\em Algebraic \& Geometric Topology}, 17(2):869--893, 2017.

\bibitem{kochen}
Simon Kochen and E.~P. Specker.
\newblock The problem of hidden variables in quantum mechanics.
\newblock {\em Journal of Mathematics and Mechanics}, 17(1):59--87, 1967.

\bibitem{bell1966problem}
John~S Bell.
\newblock On the problem of hidden variables in quantum mechanics.
\newblock {\em Reviews of Modern Physics}, 38(3):447, 1966.

\bibitem{cleve2014characterization}
Richard Cleve and Rajat Mittal.
\newblock Characterization of binary constraint system games.
\newblock In {\em International Colloquium on Automata, Languages, and
  Programming}, pages 320--331. Springer, 2014.

\bibitem{ACT12}
A.~Adem, F.~R. Cohen, and E.~Torres~Giese.
\newblock Commuting elements, simplicial spaces and filtrations of classifying
  spaces.
\newblock {\em Math. Proc. Cambridge Philos. Soc.}, 152(1):91--114, 2012.

\bibitem{Coho}
Cihan Okay, Sam Roberts, Stephen~D Bartlett, and Robert Raussendorf.
\newblock Topological proofs of contextuality in qunatum mechanics.
\newblock {\em Quantum Information \& Computation}, 17(13-14):1135--1166, 2017.

\bibitem{okay2020homotopical}
Cihan Okay and Robert Raussendorf.
\newblock Homotopical approach to quantum contextuality.
\newblock {\em Quantum}, 4:217, 2020.

\bibitem{O16}
Cihan Okay.
\newblock Spherical posets from commuting elements.
\newblock {\em J. Group Theory}, 21(4):593--628, 2018.

\bibitem{AG15}
A.~Adem and J.~M. G\'omez.
\newblock A classifying space for commutativity in {L}ie groups.
\newblock {\em Algebr. Geom. Topol.}, 15(1):493--535, 2015.

\bibitem{cohen2016survey}
Frederick~R Cohen and Mentor Stafa.
\newblock A survey on spaces of homomorphisms to Lie groups.
\newblock In {\em Configuration spaces}, pages 361--379. Springer, 2016.

\bibitem{AGV17}
O~Antol\'{i}n-Camarena, S.~Gritschacher, and B.~Villarreal.
\newblock Classifying spaces for commutativity of low-dimensional lie groups.
\newblock In {\em Mathematical Proceedings of the Cambridge Philosophical
  Society}, volume 169, pages 433--478. Cambridge University Press, 2020.

\bibitem{OW18}
Cihan Okay and Ben Williams.
\newblock On the mod-$\ell$ homology of the classifying space for
  commutativity.
\newblock {\em Algebraic \& Geometric Topology}, 20(2):883--923, 2020.

\bibitem{DV18}
Daniel~A Ramras and Bernardo Villarreal.
\newblock Commutative cocycles and stable bundles over surfaces.
\newblock In {\em Forum Mathematicum}, volume~31, pages 1395--1415. De Gruyter,
  2019.

\bibitem{ramras2018homological}
Daniel~A Ramras and Mentor Stafa.
\newblock Homological stability for spaces of commuting elements in {L}ie
  groups.
\newblock {\em International Mathematics Research Notices}, 2021(5):3927--4002,
  2021.

\bibitem{okay2020commutative}
Cihan Okay and P{\'a}l Zs{\'a}mboki.
\newblock Commutative simplicial bundles.
\newblock {\em arXiv preprint arXiv:2001.04052}, 2020.

\bibitem{gritschacher2019commuting}
Simon Gritschacher and Markus Hausmann.
\newblock Commuting matrices and Atiyah's real {$K$}-theory.
\newblock {\em Journal of Topology}, 12(3):832--853, 2019.

\bibitem{Mermin}
N~David Mermin.
\newblock Hidden variables and the two theorems of {J}ohn {B}ell.
\newblock {\em Reviews of Modern Physics}, 65(3):803, 1993.

\bibitem{Gri17}
Simon Gritschacher.
\newblock The spectrum for commutative complex {$K$}-theory.
\newblock {\em Algebraic \& Geometric Topology}, 18(2):1205--1249, 2018.

\bibitem{adams1974stable}
John~Frank Adams.
\newblock  Stable homotopy and generalised homology.
\newblock {\em University of Chicago press}, 1974.

\bibitem{switzer2017algebraic}
Robert~M. Switzer.
\newblock  Algebraic topology---homotopy and homology.
\newblock {\em Classics in Mathematics. Springer-Verlag}, Berlin, 2002.
\newblock Reprint of the 1975 original.

\bibitem{beaudry2018guide}
Agnes Beaudry and Jonathan~A Campbell.
\newblock A guide for computing stable homotopy groups.
\newblock {\em Topology and quantum theory in interaction}, 718:89--136, 2018.

\bibitem{adams1978infinite}
John~Frank Adams.
\newblock  Infinite loop spaces.
\newblock Number~90. {\em Princeton University Press}, 1978.

\bibitem{segal1974categories}
Graeme Segal.
\newblock Categories and cohomology theories.
\newblock {\em Topology}, 13(3):293--312, 1974.

\bibitem{bousfield1978homotopy}
Aldridge~K Bousfield and Eric~M Friedlander.
\newblock Homotopy theory of {$\Gamma$}-spaces, spectra, and bisimplicial sets.
\newblock In {\em Geometric applications of homotopy theory II}, pages 80--130.
  Springer, 1978.

\bibitem{schwede2018global}
Stefan Schwede.
\newblock  Global homotopy theory, volume~34.
\newblock {\em Cambridge University Press}, 2018.

\bibitem{schwede1999stable}
Stefan Schwede.
\newblock Stable homotopical algebra and $\Gamma$-spaces.
\newblock In {\em Mathematical Proceedings of the Cambridge Philosophical
  Society}, volume 126, pages 329--356, 1999.

\bibitem{bruner2003connective}
Robert~Ray Bruner and John Patrick~Campbell Greenlees.
\newblock   The connective {$K$}-theory of finite groups, volume 165.
\newblock {\em American Mathematical Soc.}, 2003.

\bibitem{hashimoto1983connective}
Shin Hashimoto.
\newblock On the connective {$K$}-homology groups of the classifying spaces
  {$B\mathbb{Z}/p^r$}.
\newblock {\em Publications of the Research Institute for Mathematical
  Sciences}, 19(2):765--771, 1983.

\bibitem{lima1960}
Elon~L. Lima.
\newblock Stable {P}ostnikov invariants and their duals.
\newblock {\em Summa Brasil. Math.}, 4:193--251, 1960.

\bibitem{casacuberta2005homotopical}
Carles Casacuberta and Javier Guti{\'e}rrez.
\newblock Homotopical localizations of module spectra.
\newblock {\em Transactions of the American Mathematical Society},
  357(7):2753--2770, 2005.


\bibitem{adem2002topics}
Alejandro Adem and James~F Davis.
\newblock Topics in transformation groups.
\newblock {\em Handbook of geometric topology}, pages 1--54, 2002.

\bibitem{slofstra2019tsirelson}
William Slofstra.
\newblock Tsirelson’s problem and an embedding theorem for groups arising
  from non-local games.
\newblock {\em Journal of the American Mathematical Society}, 2019.

\bibitem{ji2020mip}
Zhengfeng Ji, Anand Natarajan, Thomas Vidick, John Wright, and Henry Yuen.
\newblock {MIP*= RE}.
\newblock {\em arXiv preprint arXiv:2001.04383}, 2020.

\bibitem{mermin1990simple}
N~David Mermin.
\newblock Simple unified form for the major no-hidden-variables theorems.
\newblock {\em Physical review letters}, 65(27):3373, 1990.

\bibitem{peres1990incompatible}
Asher Peres.
\newblock Incompatible results of quantum measurements.
\newblock {\em Physics Letters A}, 151(3-4):107--108, 1990.

\bibitem{kitaev2013classification}
Alexei Kitaev.
\newblock On the classification of short-range entangled states. {\em {T}alk at
  {S}imons {C}enter}, 2013.

\bibitem{marcolli2019gamma}
Matilde Marcolli.
\newblock Gamma spaces and information.
\newblock {\em Journal of Geometry and Physics}, 140:26--55, 2019.

\bibitem{Slo}
Richard Cleve, Li~Liu, and William Slofstra.
\newblock Perfect commuting-operator strategies for linear system games.
\newblock {\em Journal of Mathematical Physics}, 58(1):012202, 2017.

\bibitem{qassim2020classical}
Hammam Qassim and Joel~J Wallman.
\newblock Classical vs quantum satisfiability in linear constraint systems
  modulo an integer.
\newblock {\em Journal of Physics A: Mathematical and Theoretical},
  53(38):385304, 2020.

\bibitem{TwoDim}
William~A Bogley, SJ~Pride, C~Hog-Angeloni, W~Metzler, and AJ~Sieradski.
\newblock Two-dimensional homotopy and combinatorial group theory.
\newblock {\em London Mathematical Society Lecture Note Series}, 197, 1993.

\bibitem{mccleary2001user}
John McCleary.
\newblock   A user's guide to spectral sequences.
\newblock Number~58. {\em Cambridge University Press}, 2001.


\bibitem{whitehead1949combinatorial}
John~HC Whitehead.
\newblock Combinatorial homotopy. {I}.
\newblock {\em Bulletin of the American Mathematical Society}, 55(3):213--245,
  1949.

\bibitem{spanier1989algebraic}
Edwin~H Spanier.
\newblock  Algebraic topology, volume~55.
\newblock  {\em Springer Science \& Business Media}, 1989.

\bibitem{okay2019classifying}
Cihan Okay and Daniel Sheinbaum.
\newblock Classifying space for quantum contextuality.
\newblock In {\em Annales Henri Poincar{\'e}}, volume~22, pages 529--562.
  Springer, 2021.

\bibitem{anderson1967structure}
Donald~W Anderson, Edgar~H Brown, and Franklin~P Peterson.
\newblock The structure of the spin cobordism ring.
\newblock {\em Annals of Mathematics}, pages 271--298, 1967.

\bibitem{campbell2017homotopy}
Jonathan~A Campbell.
\newblock Homotopy theoretic classification of symmetry protected phases.
\newblock {\em arXiv preprint arXiv:1708.04264}, 2017.

\bibitem{gu2014symmetry}
Zheng-Cheng Gu and Xiao-Gang Wen.
\newblock Symmetry-protected topological orders for interacting fermions:
  Fermionic topological nonlinear $\sigma$ models and a special group
  supercohomology theory.
\newblock {\em Physical Review B}, 90(11):115141, 2014.

\bibitem{kapustin2015fermionic}
Anton Kapustin, Ryan Thorngren, Alex Turzillo, and Zitao Wang.
\newblock Fermionic symmetry protected topological phases and cobordisms.
\newblock {\em Journal of High Energy Physics}, 2015(12):1--21, 2015.

\bibitem{higuera2014spaces}
Galo Higuera~Rojo.
\newblock Spaces of homomorphisms and commuting orthogonal matrices.
\newblock {\em  PhD thesis}, University of British Columbia, 2014.

\bibitem{bruner2010connective}
Robert~Ray Bruner and John Patrick~Campbell Greenlees.
\newblock  Connective real $ K $-theory of finite groups.
\newblock Number 169. {\em American Mathematical Soc.}, 2010.

\bibitem{ravenel2003complex}
Douglas~C Ravenel.
\newblock Complex cobordism and stable homotopy groups of spheres.
\newblock {\em American Mathematical Soc.}, 2003.

\bibitem{liulevicius1963theorem}
Arunas Liulevicius.
\newblock A theorem in homological algebra and stable homotopy of projective
  spaces.
\newblock {\em Transactions of the American Mathematical Society},
  109(3):540--552, 1963.

\end{thebibliography}
\end{document}